\documentclass[preview, 3p, 12pt]{elsarticle}
\usepackage{soul}
\usepackage{graphicx, graphics}
\usepackage{syntonly}
\usepackage{amssymb, amsmath, amsfonts,nccmath}
\usepackage{subfigure}
\usepackage{multirow}
\usepackage{color}
\usepackage{hyperref}

\usepackage{amsthm}

\usepackage[T1]{fontenc}
\usepackage{txfonts}
\usepackage{bm}
\usepackage[english]{babel}
\usepackage{lineno}

\providecommand{\etal}{ et al }

\newcommand{\sech}{\mbox{sech}}
\newcommand{\erf}{\mbox{erf}}
\newcommand{\ud}{\mbox{d}}

\newtheorem{theorem}{Theorem}

\begin{document}

\title{CDF-Generated Damage Laws: Admissibility, $\Gamma$-Convergence to Griffith Fracture, and Well-Posedness}
\author[TU]{Huilong Ren\corref{cor1}}
\ead{hlren@tongji.edu.cn}
\address[TU]{Department of Geotechnical Engineering,\\ College of Civil Engineering,Tongji University, Shanghai 200092, China}

\cortext[cor1]{Department of Geotechnical Engineering, College of Civil Engineering, Tongji University, Shanghai 200092, China.}
\begin{abstract}
We formulate a family of scalar softening laws by setting the stored-energy density $\psi(\eta)=\int_{0}^{\eta}[1-F(s)]\ud s$, where $F$ ranges over exponential, Cauchy, logistic, half-normal, Gudermannian, hypergeometric, radical, rational, piece-wise, and rapid-decay cumulative-distribution functions (CDFs). We prove that every such law yields a degradation map that is monotone, bounded, and dissipative, rendering the associated hyperelastic material thermodynamically admissible. Working directly in spatial dimensions $d=2,3$, we establish compactness and $\Gamma$-convergence of the CDF-based energies to a sharp-interface Griffith functional. We further show the existence of rate-independent quasi-static evolutions by constructing global energetic solutions that satisfy both stability and energy balance. These analytical results provide a rigorous bridge between the probabilistic damage formulation and Griffith–type fracture mechanics. One illustrative example is presented to show the effectiveness of the current damage laws. 

\end{abstract}
\maketitle

\section{Introduction}
Continuum Damage Mechanics (CDM) provides a rigorous framework for modeling progressive material degradation in solids. In classical CDM, a scalar damage variable $d$ is introduced to quantify the loss of stiffness from 0 (undamaged) to 1 (fully broken). This variable enters the constitutive law multiplicatively, for example via $\bm \sigma=(1-d)\mathbb D:\bm \varepsilon$, to reduce the effective stress and stiffness as damage grows. Such formulations, pioneered by Kachanov \cite{kachanov1958rupture,kachanov1986introduction}, Chaboche\cite{chaboche1988continuum}, Lemaitre \cite{lemaitre2012course}, and others \cite{murakami2012continuum,shojaei2014multi,lemaitre2006engineering,ochsner2016continuum,murakami1988mechanical,zhang2010continuum}, ensure that material softening is captured in a continuum setting. A key requirement in any damage model is irreversibility: damage can only accumulate but not heal under increasing loads. Correspondingly, CDM models typically introduce an internal variable (or history parameter) with Kuhn-Tucker loading-unloading conditions to enforce monotonic growth of $d$. These conditions, while effective, add complexity by requiring evolution equations and tracking of peak values to maintain thermodynamic consistency (i.e. non-negative energy dissipation). There is a clear incentive to seek simpler constitutive mappings that inherently satisfy the necessary properties of damage evolution - monotonicity, saturation (boundedness), and smoothness - without introducing ad-hoc internal variables.

As a fundamental mathematical framework \cite{wilce2002quantum}, probability theory has far-reaching implications across various disciplines.
In particular, probabilistic methods have proven effective when dealing with materials exhibiting inherent heterogeneity, such as natural and composite materials \cite{bazant1991statistical,mulla2022fluctuation,guy2012probabilistic}. By leveraging principles of probability theory, researchers can develop statistical models that capture the distribution of flaws or critical defects within these complex systems. Notably, Guy et al. have employed Weibull-based models to describe the initiation probabilities associated with crack formation in materials \cite{guy2012probabilistic}, thereby establishing a direct link between material properties and fracture behavior. Furthermore, probability theory has been instrumental in developing statistical models of fatigue crack growth as a function of stress amplitude \cite{suresh1998fatigue}. 

While Weibull-based models have been widely used to describe material strength distributions \cite{weibull1951statistical,bazant2019fracture} and statistical formulations of damage evolution were explored in works by Krajcinovic \cite{krajcinovic1996damage}, these models typically require postulated probability laws and are not tightly coupled with continuum variational principles. In contrast, our proposed framework derives damage models directly from probability distributions in a thermodynamically consistent form.

In this work, we propose a unified and physically interpretable class of damage models derived from cumulative distribution functions (CDFs) of the strain energy density. Specifically, the stored energy is defined by
\begin{align}
\psi(\eta) = \int_0^\eta [1 - F(s)]\, \mathrm{d}s,
\end{align}
where \( \eta = \phi(\nabla u) \) is the elastic energy density and \( F(s) \) is a cumulative distribution function. This formulation admits a wide range of classical and non-classical softening laws, including exponential, Cauchy, logistic, half-normal, Gudermannian, rational, and rapidly decaying types. Each distribution reflects different statistical assumptions on the strength and distribution of microstructural defects.

This CDF-based approach has several advantages. First, it provides a clear statistical interpretation: \( F(s) \) represents the probability that a material point has failed under energy density \( s \). The energy degradation is thus derived by integrating the survival function \( 1 - F(s) \), resulting in a damage evolution law with a strong physical basis. Second, the resulting degradation functions are automatically bounded, monotone, and strictly dissipative, ensuring thermodynamic consistency of the constitutive model.

The main objective of this paper is to establish a rigorous mathematical foundation for this class of CDF-generated damage models. Our contributions are threefold:

\begin{itemize}
  \item We prove that each CDF-based softening law produces a thermodynamically admissible energy functional in the sense of strict energy dissipation, convexity in strain, and bounded degradation.
  \item We analyze the variational structure of the model and establish compactness and \(\Gamma\)-convergence of the energy functionals to the Griffith fracture energy in spatial dimensions \( \mathbb{R}^d \), for \( d = 2,3 \). This provides a precise link between the smooth damage model and classical sharp-interface fracture.
  \item We demonstrate the existence of quasi-static evolutions under a rate-independent energetic formulation. The resulting solution concept is robust and consistent with modern variational fracture theories.
\end{itemize}

To illustrate the applicability of the proposed model, we present a single numerical benchmark in two dimensions using five representative CDF-based damage laws. We emphasize that a complete finite element implementation and systematic computational validation are presented in a separate companion paper \cite{RenCDFfem2025}, where convergence behavior, mesh sensitivity, and comparison with phase-field models are analyzed.

Overall, this work provides a statistical-to-continuum bridge between flaw population models and damage mechanics, establishing a flexible and mathematically sound foundation for computational fracture. It opens avenues for model calibration based on probabilistic microstructure data and motivates further investigation into the thermodynamic structure of damage evolution. 

The remainder of the paper is organized as follows. Section 2 reviews standard probability distributions and introduces new CDFs tailored to constitutive modeling. Section 3 develops the general CDF-based damage framework, and Section 4 formulates and analyzes several specific models. Section 5 establishes the $\Gamma$-convergence to Griffith fracture, while Section 6 proves existence of rate-independent quasi-static evolutions of the CDFs-based damage laws. Section 7 presents a 2D numerical example, and Section 8 concludes with a summary and outlook. Together, the proposed framework offers a probabilistically grounded, variationally consistent approach to damage modeling.

\section{Some new probability distributions}
In probability theory, the fundamental concepts are the probability density function and cumulative distribution function. For continuous probability distribution, PDF is a non-negative function and CDF is a monotonically non-decreasing, right-continuous function, satisfying $ \lim _{x\to -\infty }F(x)=0\,$ and $ \lim _{x\to+\infty }F(x)=1$ for a random variable $x$ defined in ($-\infty,+\infty$) \cite{gnedenko2018theory}.
The types of probability distributions in probability theory are abundant, and each of them represents a unique distribution. For simplicity, we take the exponential distribution for an example. 

The exponential distribution models the time between independent events that occur at a constant average rate, with its probability density function decreasing exponentially \cite{walck1996hand}. It is commonly used to describe waiting times in a Poisson process, where the memoryless property is a key characteristic.
The probability density function of the exponential distribution is expressed as:
\begin{align}
f(x; \lambda) = \lambda e^{-\lambda x}, \quad x \geq 0,
\end{align}
where $\lambda > 0$ is the rate parameter, representing the average number of occurrences per unit time.
The cumulative distribution function, which provides the probability that a random variable $X$ is less than or equal to a given value $x$, is given by:
\begin{align}
F(x; \lambda) = 1 - e^{-\lambda x}, \quad x \geq 0. \label{eq:Fexp0}
\end{align}
For a more detailed description of the exponential distribution, the reader is refered to Ref \cite{walck1996hand}. Besides the exponential distribution, other traditional probability distributions suitable for damage models include Cauchy distribution, Logistic distribution, Half-normal distribution, Chi-square distribution, which are supplemented in \ref{sec:probdata}. 

Traditional probability models are often expressed in exponential or transcendental functions. Based on the principle of probability distribution, we propose some new probability distributions based on the radical function, the piecewise function, rapid decay function, the Gudermannian function, hypergeometric function, and rational function.

\subsection{Radical distribution}
\begin{figure}[htp]
\centering
\subfigure[]{
\includegraphics[width=0.46\textwidth]{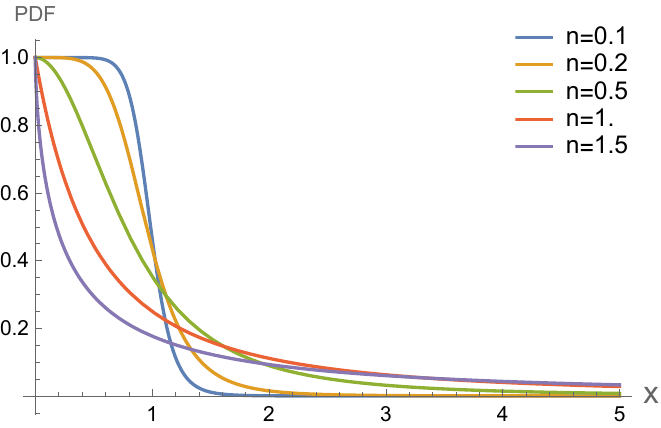}
\label{fig:raydisa}}
\subfigure[]{
\includegraphics[width=0.46\textwidth]{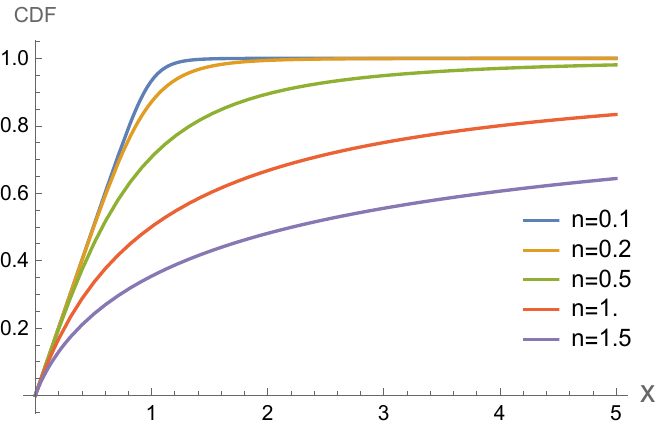}
\label{fig:raydisb}}
\caption{Radical distribution: (a) probability density function (b) cumulative distribution function}
\label{fig:Rational2}
\end{figure}

According to the requirements of CDF, we employ the radical function (or root form function) to design a new CDF of random variable $x\in [0,+\infty)$ defined as
\begin{align}
F(x; n)=\frac{x}{(1+x^{1/n})^n},\, x\geq 0 \label{eq:newProb1}
\end{align}
where $n>0$ is the scale parameter of the distribution. This CDF contains $n$-fold square root of a polynomial. It can be verified easily that $F(0;n)=0, \,\lim_{x\to +\infty} F(x; n)=1$. $F(x;n)$ is a monotonically increasing function and satisfies the requirements of a probability function. We call this new probability distribution the radical distribution or root form distribution. The probability density function is
\begin{align}
f(x; n)=F'(x;n)=\frac{1}{(1+x^{1/n})^{n+1}},\,x\geq 0,n>0.\label{eq:pdf1}
\end{align}
The PDF and CDF of this radical distribution are depicted in Figure \ref{fig:Rational2}.
The $m$-order plain moment exists if $m n<1$:
\begin{align}
\mathbb E(x^m)=\int_{0}^{+\infty} f(x;n) x^m \ud x =\frac{\Gamma(1-mn) \Gamma(n+mn)}{\Gamma(n)} \mbox{, if } m n<1
\end{align}
where $\mathbb E(x^m)$ denotes the $m$-th plain moment of the random variable.
This distribution contains many special cases, for example,
\begin{align}
&F(x, 1)=\frac{x}{1+x},\,\, F(x,2)=\frac{x}{(1+\sqrt{x})^2}, \,\,F(x; \tfrac 12)=\frac{x}{\sqrt{1+x^{2}}}\notag\\
&f(x,1)=\frac{1}{(1+x)^{2}},\,\, f(x,2)=\frac{1}{(1+\sqrt{x})^{3}}, \,\,f(x,\tfrac 12)=\frac{1}{(1+x^{2})^{3/2}}\notag
\end{align}
Using mathematical software, it is easy to verify that the PDF given by Eq.\ref{eq:pdf1} satisfies
\begin{align}
\lim_{n\to 0^+} f(x;n)=1 ,\,\,\forall x\in (0,1),\quad
\lim_{n\to 0^+} f(x;n)=0 ,\,\,\forall x>1\notag
\end{align}
Therefore, in the limit of $n\to 0$, the radical distribution degenerates to the uniform distribution in [0,1]. For this reason, the radical distribution is a generalization of the uniform distribution.

\subsection{Piecewise distribution}
\begin{figure}[htp]
\centering
\subfigure[]{
\includegraphics[width=0.46\textwidth]{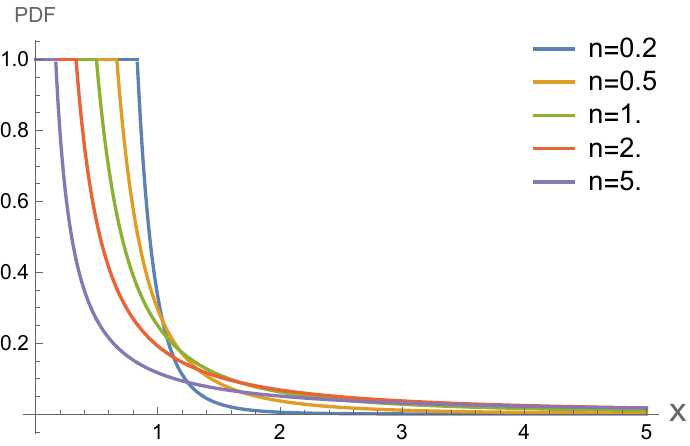}
\label{fig:raydisa}}
\subfigure[]{
\includegraphics[width=0.46\textwidth]{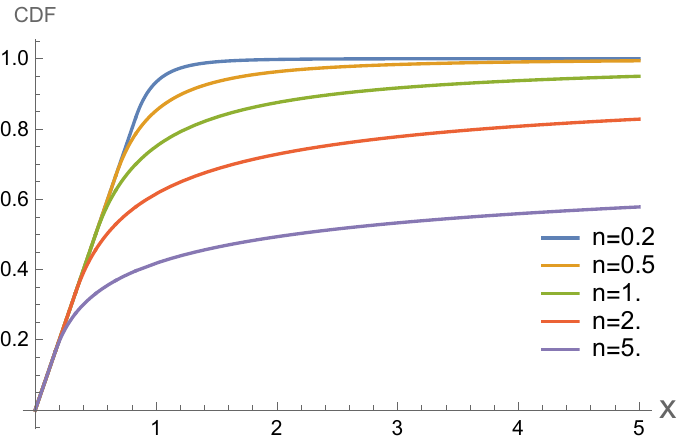}
\label{fig:raydisb}}
\caption{Piecewise distribution: (a) probability density function (b) cumulative distribution function}
\label{fig:Piecewise2}
\end{figure}

By employing a piecewise continuous function, we propose a new accumulated distribution function defined as
\begin{align}
F(x;n)=\mbox{If}\big(x\leq \frac{1}{n+1}, x, 1-\frac{{n}/(n+1)}{((n+1)x)^{1/n}}\big), \label{eq:newProb2}
\end{align}
where ``{If}(\textit{condition}, $T$, $F$)'' gives $T$ if \textit{condition} evaluates to \textit{True}, and $F$ if it evaluates to \textit{False}. It can be easily verified that $F(0;n)=0$, $\lim_{x\to +\infty} F(x;n)=1$ and $F(x;n)$ is a non-decreasing continuous function. We name this distribution as the piecewise distribution.
For this function, the probability density function is
\begin{align}
f(x;n)=F'(x;n)=\mbox{If}\big(x\leq \frac{1}{n+1}, 1, \big((n+1)x\big)^{-\frac{n+1}{n}}\big).\label{eq:pdf2}
\end{align}
The PDF and CDF of this piecewise distribution are depicted in Figure \ref{fig:Piecewise2}, which shows that the distributions become wider as $n$ increases.
The $m$-th plain moment is
\begin{align}
\mathbb E(x^m)=\int_{0}^{+\infty} f(x;n) x^m \ud x =\frac{(n+1)^{-m}}{(m+1)(1-m n)} \mbox{, if } m n<1
\end{align}
where $\mathbb E(x^m)$ denotes the $m$-th plain moment of the random variable.
This distribution contains many special cases, for example
\begin{align}
f(x,1)=\mbox{If}\big(x\leq \frac{1}{2}, 1, \frac{1}{4 x^2}\big), f(x,2)=\mbox{If}\big(x\leq \frac{1}{3}, 1, \frac{1}{\sqrt{27 x^3}}\big), f(x,\frac 12)=\mbox{If}\big(x\leq \frac{2}{3}, 1, \frac{8}{27 x^3}\big).\notag
\end{align}
Using mathematical software, it is easy to verify that the PDF given by Eq.\ref{eq:pdf2} satisfies
\begin{align}
\lim_{n\to 0^+} f(x;n)=1 ,\,\,\forall x\in (0,1);\quad
\lim_{n\to 0^+} f(x;n)=0 ,\,\,\forall x>1.\notag
\end{align}
Therefore, in the limit of $n\to 0$, the piecewise distribution degenerates to the uniform distribution in $[0,1]$. For this reason, the piecewise distribution is another generalization of continuous uniform distribution.

In this section, we introduce two novel probability distributions derived from the conditions of cumulative distribution functions in probability theory. These new distributions showcase the potential for further discoveries within the vast expanse of mathematical theory. Besides above two new distributions, we also propose four new probability distributions: rational distribution, Gudermannian distribution, continuous hypergeometric distribution, and rapid decay distribution, which are presented in \ref{sec:probdatanew}. In the subsequent sections, we will demonstrate the potent applications of these new distributions in the realm of damage mechanics.

\section{Mathematical framework}

\subsection{Small-strain kinematics and elastic energy}
We consider an elastic solid in the small-strain regime. Let $\bm {u}(\bm {x})$ be the displacement field defined on the body domain $\Omega$. The infinitesimal strain tensor is given by the symmetric gradient of the displacement:
\begin{align}
\bm {\varepsilon}(\bm {u})=\nabla^s \bm {u}=\frac{1}{2}\left(\nabla \bm {u}+(\nabla \bm {u})^T\right),
\end{align}
which assumes small deformations so that geometric nonlinearities are neglected. Under linear isotropic elasticity, the undamaged (virgin) material has a Hookean elastic energy density
\begin{align}
\phi _0(\bm {\varepsilon})=\frac{1}{2} \bm {\varepsilon}: \mathbb{D}: \bm {\varepsilon}, 
\end{align}
where $\mathbb{D}$ is the fourth-order elastic stiffness tensor (for isotropic materials, $\mathbb{D}$ can be expressed in terms of the Lamé constants or Young’s modulus $E$ and Poisson’s ratio $\nu$). This $\phi _0$ represents the stored energy per unit volume in the absence of damage.
\subsection{Variational principle of total energy}
The equilibrium configuration of the body is determined by the principle of minimum total potential energy. We introduce the total energy functional $E[\bm {u}]$, which is the difference between the internal (strain) energy and the work of external forces. In the presence of damage, the internal energy density will be modified (as detailed below), but the governing principle remains:
\begin{align}
E[\bm {u}]=\int_{\Omega} \phi(\bm {\varepsilon}(\bm {u})) \ud V-\int_{\Omega} \bm {b} \cdot \bm {u} \ud V-\int_{\Gamma_t} \overline{\bm {t}} \cdot \bm {u} \ud S,\label{eq:PiEnergy}
\end{align}
where $\phi(\bm {\varepsilon})$ is the damaged elastic energy density (to be defined), $\bm {b}$ is a body force per unit volume, and $\bar{\bm {t}}$ is a prescribed traction on the Neumann boundary $\Gamma_t$. The principle of stationary energy $\delta E=0$ will yield the equilibrium (Euler-Lagrange) equations for the displacement field.

\subsection{Spectral decomposition into tensile and compressive parts}
A key feature of the damage model is the treatment of tension and compression differently. To avoid unphysical damage growth under compression (crack closure), we split the strain energy into tensile (positive) and compressive (negative) contributions. This is achieved via a spectral decomposition of the strain tensor. Let ${\varepsilon_i, \bm {n}_i}$ for $i=1,2,3$ denote the principal strains (eigenvalues of $\bm {\varepsilon}$) and their orthonormal principal directions. We decompose $\bm {\varepsilon}$ as
\begin{align}
\bm {\varepsilon}=\sum_{i=1}^3 \varepsilon_i \bm {n}_i \otimes \bm {n}_i=\bm {\varepsilon}^{+}+\bm {\varepsilon}^{-},
\end{align}
where the tensile part $\bm {\varepsilon}^+$ contains only the positive principal strains and the compressive part $\bm {\varepsilon}^-$ contains only the negative principal strains. In other words, $\bm {\varepsilon}^+ = \sum_{i} \langle\varepsilon_i\rangle_+\bm {n}_i\otimes \bm {n}_i$ and $\bm {\varepsilon}^- = \sum_{i} \langle\varepsilon_i\rangle_-\bm {n}_i\otimes \bm {n}_i$, where $\langle x\rangle_\pm=(x\pm |x|)/2$. Correspondingly, we define the tensile elastic energy and compressive elastic energy as
\begin{align}
\phi ^{+}(\bm {\varepsilon})=\frac{1}{2} \bm {\varepsilon}^{+}: \mathbb{D}: \bm {\varepsilon}^{+}, \quad \phi ^{-}(\bm {\varepsilon})=\frac{1}{2} \bm {\varepsilon}^{-}: \mathbb{D}: \bm {\varepsilon}^{-} . 
\end{align}

Thus $\phi _0(\bm {\varepsilon}) = \phi ^+(\bm {\varepsilon}) + \phi ^-(\bm {\varepsilon})$. Only the positive part $\phi ^+$ will be subject to damage, reflecting the physical fact that micro-crack growth is driven by tensile strain energy, while compressive states do not lead to further damage

\subsection{Admissible damage mappings: A sufficiency theorem}
We now formalise the requirements a scalar degradation map must satisfy and prove that any CDF is admissible.
\begin{theorem}
 (Admissible degradation map).
Let $g:[0,\infty)\to [0,1]$ be $C^1$. Suppose
\begin{align}
\text {   (i) } g(0)=0, \quad \text {   (ii) } \lim _{\eta \rightarrow \infty} g(\eta)=1, \quad \text {   (iii) } g^{\prime}(\eta)\ge 0\quad \forall \eta>0 \label{eq:g16}
\end{align}
Then $d:=g(\eta)$ is (a) bounded, (b) monotonically non-decreasing, and (c) yields a finite fracture energy
\begin{align}
\int_0^{\infty}(1-g(\eta)) \,\ud \eta<\infty. \label{eq:g17}
\end{align}
Conversely, any $g$ that violates Eq.\ref{eq:g16} fails to satisfy at least one of the properties (a)–(c).
\end{theorem}

\noindent\textbf{Proof.} (a)-(b) follow directly from Eq.\ref{eq:g16}. For (c)  boundedness plus monotone convergence implies $(1-g(\eta))$ is integrable on $[0,\infty)$.

\textbf{Remark 1}. A cumulative distribution function satisfies Eq.\ref{eq:g16} by definition, hence every CDF is an admissible degradation map. Conversely, an arbitrary polynomial or cubic-spline softening law may violate monotonicity or boundedness and is therefore not automatically admissible. This theorem justifies our exclusive use of CDF-based damage laws.

\subsection{CDF-based damage law for tensile energy}

To model progressive material degradation, we introduce a damage formulation that acts on the tensile energy component $\phi ^+$. The damage is characterized by a continuous variable $d$ that evolves from 0 (undamaged material) to 1 (fully broken material). Instead of prescribing an arbitrary evolution law, we derive $d$ from a CDF type mapping of the tensile energy. Specifically, we choose a monotonically increasing scalar function $F:[0,\infty)\to[0,1]$ with $F(0)=0$ and $F(\infty)=1$ - having the characteristic S-shape of a CDF. The damage variable at a material point is then defined as
\begin{align}
d=F\left(\phi ^{+}(\bm {\varepsilon})\right),
\end{align}
so that higher tensile energy density produces more damage. In the initial elastic regime $\phi ^+\approx 0$, we have $d\approx 0$, meaning the material is essentially intact. As $\phi ^+$ grows, $d$ increases toward 1 in a smooth sigmoidal manner. The precise form of $F$ can be chosen based on experimental calibration; typical choices include logistic-type functions or error-function (Gaussian CDF) shapes, which ensure a gradual transition from $d\approx 0$ to $d\approx 1$. This CDF-based approach guarantees that $d$ remains bounded between 0 and 1 and evolves in a stable, gradual fashion rather than instantaneously. The damage influence is introduced into the energy functional by reducing the tensile energy capacity. We define an effective tensile energy density $\psi(\phi ^+)$ as the portion of tensile energy that remains stored in the material after damage. Consistent with the definition of $d$, we take
\begin{align}
\boxed{\psi \left(\phi ^{+}\right)=\int_0^{\phi ^{+}}(1-F(s)) \ud s}, \label{eq:prob8}
\end{align}
so that the fraction $1-d = 1-F(\phi ^+)$ plays the role of a stiffness reduction factor. Indeed, differentiating this relation yields ${\ud{\psi  }}/{\ud \phi ^+} = 1 - F(\phi ^+) = 1-d$, indicating that an infinitesimal increase $\ud \phi ^+$ in tensile energy leads to an increase $\ud{\psi  }= (1-d) \ud \phi ^+$ in stored energy. In the undamaged limit ($d=0$) we have $\psi = \phi ^+$, while as $d\to 1$ the stored tensile energy increment tends to zero (the material can no longer sustain additional tensile energy). 

In addition,  for many forms $F(s)$, Eq.\ref{eq:prob8} can be integrated as rational polynomials or other simple function expressions, which shall be presented in the subsequent sections.

The total damaged energy density is then defined as
\begin{align}
 {\phi}(\bm {\varepsilon})=\psi\left(\phi ^{+}(\bm {\varepsilon})\right)+\phi ^{-}(\bm {\varepsilon}),\label{eq:esplitpm}
\end{align}
combining the degraded tensile part and the undegraded compressive part. By construction, $\phi(\bm {\varepsilon})$ is a concave function of $\phi ^+$ that saturates as $\phi ^+\to\infty$, reflecting the approach to complete material failure. The function $F$ (and its parameters) can be tuned such that the area under the stress-strain curve equals the material’s critical fracture energy $G$, thereby ensuring consistency with the material’s toughness.

\subsection{Effective stress and equilibrium equations}
We now derive the governing equations by taking the first variation of the total energy functional. For a virtual displacement $\delta\bm {u}$, the variation $\delta E$ is
\begin{align}
\delta E=\int_{\Omega}\bm \sigma: \delta \bm {\varepsilon} \ud V-\int_{\Omega} \bm {b} \cdot \delta \bm {u} \ud V-\int_{\Gamma_t} \overline{\bm {t}} \cdot \delta \bm {u} \ud S,
\end{align}
where $\bm \sigma$ is the Cauchy stress tensor associated with the damaged energy density $\phi$. The stress can be obtained by differentiating $\phi$ with respect to the strain:
\begin{align}
\bm \sigma(\bm {\varepsilon})=\frac{\partial \phi}{\partial \bm {\varepsilon}}=\frac{\partial \psi}{\partial \bm {\varepsilon}}+\frac{\partial \phi ^{-}}{\partial \bm {\varepsilon}} .
\end{align}
Using the chain rule together with the previous definitions, this yields the effective stress in terms of the undamaged (elastic) stress contributions:
\begin{align}
\bm\sigma(\bm {\varepsilon})=\underbrace{\left(1-F\left(\phi ^{+}(\bm {\varepsilon})\right)\right) \sigma_0^{+}(\bm {\varepsilon})}_{\text {degraded tensile stress}}+\underbrace{\sigma_0^{-}(\bm {\varepsilon})}_{\text {compressive stress}} .
\end{align}

Here $\bm \sigma^+_0 = \mathbb{D}:\bm {\varepsilon}^+$ and $\bm \sigma^-_0 = \mathbb{D}:\bm {\varepsilon}^-$ are the stress contributions that would arise from the tensile and compressive strains if the material were undamaged. Thus, the tensile part of the stress is scaled by the factor $1-d = 1-F(\phi ^+)$, while the compressive part is unaffected by damage. This form of $\bm \sigma$ can also be written compactly as $\bm \sigma = \mathbb{D} : (\bm {\varepsilon}^- + (1-d),\bm {\varepsilon}^+)$, showing that the stiffness in tension is reduced by $(1-d)$ whereas in compression it remains full. Imposing $\delta E=0$ for all admissible $\delta\bm {u}$ and performing integration by parts on the internal energy term, we obtain the strong-form equilibrium equation:
\begin{align}
\nabla \cdot \bm\sigma(\bm {\varepsilon}(\bm u))+\bm {b}=\bm {0} \quad \text { in } \Omega, 
\end{align}
subject to the usual boundary conditions $\bm \sigma,\bm {n} = \bar{\bm {t}}$ on $\Gamma_t$ (traction boundary) and $\bm {u}=\bm {0}$ on any fixed-displacement boundary (or a prescribed $\bm {u}$ on $\Gamma_u$). The equilibrium balance can be seen more explicitly by substituting the expression for $\bm \sigma$ from above:
\begin{align}
\nabla \cdot\left((1-d(\phi^+))\bm \sigma_0^{+}(\bm {\varepsilon})+\bm\sigma_0^{-}(\bm {\varepsilon})\right)+\bm {b}=\bm {0} .
\end{align}
Here $d(\phi^+) = F(\phi ^+(\bm {\varepsilon}(x)))$ means that the damage field varies over $\Omega$. This is the Euler-Lagrange equation corresponding to the energy minimization principle. It represents a modified form of the linear momentum balance, in which the material’s stress-carrying capacity in tension is locally reduced in proportion to the damage $d(\phi^+)$.

\subsection{Probabilistic micro-mechanical interpretation}
Classical statistical fracture theory (Weibull 1951 \cite{weibull1951statistical}; Bažant \& Planas 1998 \cite{bazant2019fracture}) interprets damage as the volume fraction of micro-cells that have failed when a local driving quantity exceeds their random strength. Let
\begin{align}
X \sim F_X(x) \in[0,1], \quad x \geq 0 \notag
\end{align}
denote the random micro-strength of material cells, and let
\begin{align}
\eta(\mathbf{x}, t):=\phi^{+}(\boldsymbol{\varepsilon}(\mathbf{x}, t))\notag
\end{align}
be the local tensile strain-energy density that drives failure (Eq.\ref{eq:prob8}). A micro-cell fails when $\eta \ge X$. The expected fraction of failed cells at $(x,t)$ is therefore
\begin{align}
d(\mathbf{x}, t)=\operatorname{Pr}[X \leq \eta]=F_X(\eta(\mathbf{x}, t)) \label{eq:g15}
\end{align}
Hence the continuum damage variable is identified with the cumulative distribution function of the underlying strength statistics. Equation \ref{eq:g15} establishes a direct physical link between probability theory and macroscopic damage: $d$ is the mesoscale average of a Bernoulli failure indicator, and $F(x)$ is the CDF of micro-strength. The specific functional forms arise from (i) classical strength distributions (exponential, half-normal, Cauchy, logistic) and (ii) mathematically admissible generalizations tailored in Section 2.

\subsection{Thermodynamic consistency and damage irreversibility}

For completeness we verify that the CDF-based model obeys the Clausius–Duhem inequality. Using the free energy
\begin{align}
\phi(\bm\varepsilon)=\psi\left(\phi^{+}(\bm\varepsilon)\right)+\phi^{-}(\bm\varepsilon), \quad \psi^{+}(\eta):=\int_0^\eta(1-d(s)) \ud s\notag
\end{align}
and the energy split in Eq.\ref{eq:esplitpm}, the local dissipation is
\begin{align}
\mathcal{D}=\bm{\sigma}: \dot{\bm\varepsilon}-\dot{\phi}=(1-d) \bm\sigma^{+}: \dot{\bm\varepsilon}^{+}-(1-d) \bm\sigma^{+}: \dot{\bm\varepsilon}^{+}+\dot{d} \phi^{+}=\dot{d} \phi^{+}(\bm\varepsilon) \geq 0 
\end{align}
because $d'=F'_{X}(\eta)>0$ (monotone CDF) and $\phi^+\ge 0$. During unloading we enforce
\begin{align}
\eta(\mathbf{x}, t)=\max _{s \leq t} \phi^{+}(\bm\varepsilon(\mathbf{x}, s))
\end{align}
so that $\dot{d}\ge 0$ everywhere. Hence the model is thermodynamically admissible and damage evolution is irreversible.

The monotonic CDF mapping $F$ thus builds in a one-way evolution: once a material point is damaged, the reduction in stiffness persists, and any further loading will only increase the damage. This approach yields a thermodynamically consistent model of degrading elasticity, with no spurious energy generation (the second law is satisfied) and an intrinsic coupling between damage growth and energy absorption (dissipation) equal to the material’s fracture energy.

\section{From probability models to damage models}

In this section, we leverage several probability distributions, including the exponential, Cauchy, radical, and piecewise distributions, to demonstrate the transformation of probability models into innovative damage models. Based on these distributions, we propose new damage models, such as the exponential damage model, Cauchy-type damage model, among others. These models are unique in that they do not explicitly depend on a damage variable, which is instead utilized solely for descriptive purposes or post-processing. 

\subsection{Damage model based on CDF $F(x)=\min(x^n,1)$}
If we select the distribution function as $F(x)=\min(x^n,1), \,n>0$, which satisfies the principle of CDFs. We assume the expression of $\psi$ as
\begin{align}
\psi({\phi^+})=\int_0^{\phi^+} [1-F(\gamma s)]\ud s
\end{align}
where $\gamma$ is the coefficient to be determined. Set the saturation of $\max(\psi)=G/\ell$ leads to $\gamma=\frac{n\ell}{(n+1)G}$. Integrating above $\psi$ explicitly, we obtain the energy denstiy involving damage
\begin{align}
\psi({\phi^+})=\begin{cases} {\phi^+}\Big(1-\frac{1}{(n+1)}\big(\frac{n\ell\phi^+}{(n+1)G}\big)^n\Big) &\mbox{ if } \phi^+\leq \frac{(n+1)G}{n\ell}\\
G/\ell &\mbox{ otherwise }
\end{cases}\label{eq:Funxn}
\end{align}
The variation of $\psi(\phi^+)$ is
\begin{align}
\delta \psi(\phi^+)=\max\Big(0,1-\big(\frac{n\ell \phi^+}{(n+1)G}\big)^n\Big)\delta \phi^+
\end{align}

Based on the energy functional in Eq.\ref{eq:PiEnergy} and the damage energy density $\psi$ given by Eq.\ref{eq:Funxn}, the variation of $E $ yields
\begin{align}
\delta E &=\int_\Omega \delta \psi+ \delta \phi^--\bm b\cdot \delta \bm u \ud V\notag\\
&=\int_\Omega \max\Big(0,1-\big(\frac{n\ell \phi^+}{(n+1)G}\big)^n\Big) \bm \sigma^+ :\delta \bm \varepsilon+ \bm \sigma^- :\delta \bm \varepsilon-\bm b \cdot \delta \bm u \ud V
\end{align}
where $\bm \sigma^\pm=\frac{\partial \phi^\pm}{\partial \bm \varepsilon}$.
Via integration by parts, the governing equations of elastic solid with fracture is
\begin{align}
\nabla \cdot ( \max\Big(0,1-\big(\frac{n\ell \phi^+}{(n+1)G}\big)^n\Big) \bm \sigma^++\bm \sigma^-)+\bm b=\bm 0 \label{eq:cdfFxn}
\end{align}
Although the above governing equation based on CDF $F(x)=\min(x^n,1)$ has no explicit dependence on the damage variable, it incorporates the damage evolution in an elegant way. Follow the idea of damage mechanics, the damage variable for the purpose of post-processing is
\begin{align}
d=1-\max\Big(0,1-\big(\frac{n\ell \phi^+}{(n+1)G}\big)^n\Big)=\min\Big(1,\big(\frac{n\ell \phi^+}{(n+1)G}\big)^n\Big)
\end{align}

In 1D, the stress can be written as 
\begin{align}
\sigma^+=\begin{cases}
k\varepsilon\Big(1-\big(\frac{n\ell \phi^+}{(n+1)G}\big)^n\Big) & \mbox{ if }0\leq \varepsilon \leq \sqrt{\frac{2(n+1)G}{n k\ell}}  \\
0 &\mbox{ if } \varepsilon >\sqrt{\frac{2(n+1)G}{n k\ell}} 
\end{cases}
\end{align}

\begin{figure}[!htb]
\centering
\includegraphics[width=8cm]{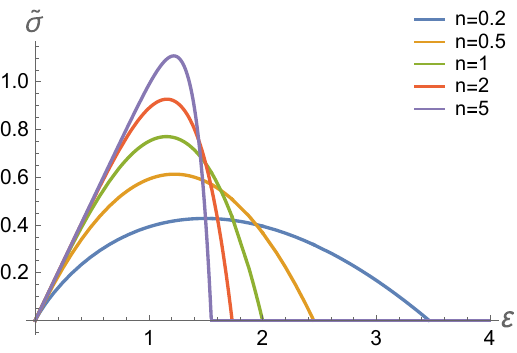}
\caption{Effective stress of damage model based on CDF $F(x)=\min(x^n,1)$ in 1D with $G/\ell=1, k=1$ and $n\in \{0.2,0.5,1,2,5\}$.}
\label{fig:FsSgeneral1}
\end{figure}

\begin{figure}[!htb]
\centering
\includegraphics[width=8cm]{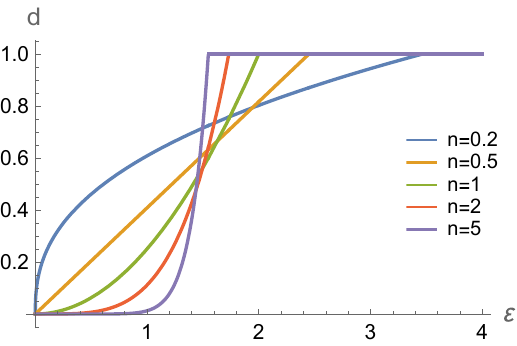}
\caption{Damage of damage model based on CDF $F(x)=\min(x^n,1)$ in 1D with $G/\ell=1, k=1$ and $n\in \{0.2,0.5,1,2,5\}$.}
\label{fig:FsDgeneral1}
\end{figure}

In the setting of 1D with assumption $G/\ell=1, k=1$ and energy density $\phi=\frac 12 k\varepsilon^2$, we plot the graphs of damage function and effective stress in Figure \ref{fig:FsDgeneral1} and Figure \ref{fig:FsSgeneral1}, respectively.
 
In Figure \ref{fig:FsSgeneral1}, it is evident that the maximal effective stress \(\tilde{\sigma}\) exhibits an increasing trend with the parameter \(n\). Specifically, for lower values of \(n\), the material experiences an early onset of stress softening, albeit at a reduced rate. Conversely, at higher \(n\) values, the material achieves a greater peak effective stress, and the softening effect is more pronounced and rapid upon reaching the maximum stress level, indicative of a brittle-like softening behavior.

Figure \ref{fig:FsDgeneral1} illustrates the progression of damage as a function of positive strain. For lower \(n\) values, the damage evolution is characterized by rapid growth during the initial stages, followed by a deceleration in the later stages, which is consistent with a more ductile damage progression. In contrast, at higher \(n\) values, the damage evolution is initially slower, but accelerates significantly in the later stages, reflecting a more brittle damage progression.

\subsection{Exponential damage model based on exponential distribution}
For the exponential distribution function $F(x)=1-\exp(-x)$, we assume the expression of $\psi$ as
\begin{align}
\psi({\phi^+})=\int_0^{\phi^+} [1-F(\gamma s)] \ud s
\end{align}
where $\gamma$ is the coefficient to be determined. According to saturation condition of $\max(\psi)=G/\ell$, we can solve $\gamma=\ell/G$. Then the damage function can be written as
\begin{align}
\psi({\phi^+})=\int_0^{\phi^+} [1-(1-\exp(-\frac{s\ell}{G})] \ud s=\frac{G}{\ell}(1-e^{-\frac{\ell}{G}\phi^+})
\end{align}
We call this new damage model the exponential damage model. 

Actually, the exponential damage model can be derived by a simple algebraic replacement as well. Based on the exponential model, we replace the random variable $x$ with $\phi$ in Eq.\ref{eq:Fexp0}, i.e. $1-e^{-\lambda x}\to \frac{G}{\ell}(1-e^{-\frac{\ell}{G}\phi})$, and define the fracture energy density as follows:
\begin{align}
\psi=\frac{G}{\ell}(1-e^{-\frac{\ell}{G}\phi})
\end{align}
The function $\psi$ is similar to the CDF, but the new limit of $\psi$ is no longer 1 but ${G}/\ell$. 

Interesting, for the exponential distribution, based on Eq.\ref{eq:prob8}, there is $\int_0^{\phi}(1-F(\lambda s)) \ud s= F(\lambda \phi)/\lambda$, which confirms above replacement. 

In the following derivation of CDFs-based damage models, both the transformation by Eq.\ref{eq:prob8} and method by algebraic replacement can be used. For simplicity, the algebraic replacement method is preferred. 

When $\phi$ is small, we have
\begin{align}
\frac{G}{\ell}(1-e^{-\frac{\ell}{G}\phi})= \phi-\frac{\ell}{2 G}\phi^2+\frac{\ell^2}{6 G^2}\phi^3+O(\phi^4),
\end{align}
which signifies that the model effectively captures linear elasticity in the absence of damage. Figure \ref{fig:LinearApproximationPhi} shows the approximation of the exponential damage model on an elastic model. When the deformation energy satisfies $\phi\ell/G \leq 0.3$, the maximal difference between damage model and elastic model is 13.6\%. This indicates that the new damage model converges to the linear elastic model when the deformation is small. 

\begin{figure}[!htb]
\centering
\includegraphics[width=8cm]{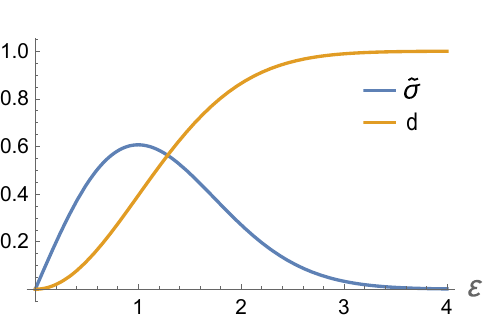}
\caption{Damage model based on exponential distribution in 1D with $G/\ell=1, k=1$: damage $d$ and effective stress $\tilde{\sigma}=e^{-\frac{\ell \phi }{G}} \sigma$.}
\label{fig:DamageCrack1}
\end{figure}

We employ the energy functional for elastic solid $E=\int_\Omega \psi+\phi^--\bm b \cdot \bm u \,\ud V$, 
where $\psi$ is the fracture energy density defined by $\psi=\frac{G}{\ell}(1-e^{-\frac{\ell}{G}\phi^+})$ and $\phi^-$ is the energy density unaffected by fracture. $\phi=\phi^++\phi^-$ is energy decomposition based on the spectral decomposition of the deformation tensor. The energy decomposition considers the fact that the tensile strength and compressive strength of a material is different.
The variation of $E$ yields
\begin{align}
\delta E&=\int_\Omega \delta \psi+ \delta \phi^--\bm b\cdot \delta \bm u \ud V\notag\\
&=\int_\Omega e^{-\frac{\ell \phi^+ }{G}} \bm \sigma^+ :\delta \bm \varepsilon+ \bm \sigma^- :\delta \bm \varepsilon-\bm b \cdot \delta \bm u \ud V
\end{align}
where $\bm \sigma^\pm=\frac{\partial \phi^\pm}{\partial \bm \varepsilon}$.
Via integration by parts, the governing equations of elastic solid with fracture is
\begin{align}
\nabla \cdot ( e^{-\frac{\ell \phi^+ }{G}} \bm \sigma^++\bm \sigma^-)+\bm b=\bm 0 \label{eq:expDamage}
\end{align}
Although the above governing equation based on exponential distribution has no explicit dependence on the damage variable, it incorporates the damage evolution in an elegant fashion. Follow the idea of damage mechanics, we introduce the damage variable for the purpose of post-processing as
\begin{align}
d=1-e^{-\frac{\ell \phi^+ }{G}}
\end{align}
For the purpose of illustration, we consider the scenario of this damage model with unit sectional area in 1 dimensional space. In the 1D setting, the elastic energy density can be simply written as $\phi=\frac 12 k \varepsilon^2$, where $k$ is the elastic modulus of the bar.
It can be easily verified that
\begin{align}
\int_{0}^{+\infty} (1-d) \sigma \,\ud \varepsilon=\frac{G}{\ell}
\end{align}
By finding the solution of $\frac{\partial \sigma}{\partial \varepsilon}=0$, the maximal stress is $\bm \sigma_{max}= \sqrt{\frac{G k}{e \ell}}$ at $\varepsilon=\sqrt{G/(\ell k)}$. The corresponding damage is $d=1-1/\sqrt{e}$.
In the setting of 1D with assumption $G/\ell=1, E=1$ and energy density $\phi=\frac 12 k\varepsilon^2$, we plot the graphs of damage function and effective stress in Figure \ref{fig:DamageCrack1}, which shows non-decreasing property of damage and the stress soften due to damage.
In previous derivation, the crack length scale $\ell$ is a parameter that remains undetermined. In order to determine the value of $\ell$, we employ the maximal stress determined from experiment. Thus, we solve $G/\ell$ from $\bm \sigma_{max}=\sqrt{\frac{k G}{e \ell}}$ as
\begin{align}
\frac{G}{\ell}=e\frac{\sigma_{max}^2}{k}. \label{eq:Gellv}
\end{align}
\subsection{Cauchy-type damage model based on Cauchy distribution}
In this damage model, we select the Cauchy distribution. Replacing the random variable $x$ in Eq.\ref{eq:Caudis} with strain energy density $\phi$ and neglecting the constant terms,
\begin{align}
{\frac {1}{\pi }}\arctan \big({\frac {x}{\gamma }}\big)+{\frac {1}{2}}\,\,\to\,\, \frac{2G}{\pi\ell} \arctan(\frac{\pi\ell}{2G}\phi),\notag
\end{align}
and we arrive at a new fracture energy density defined by
\begin{align}
\psi=\frac{2G}{\pi\ell} \arctan(\frac{\pi\ell}{2G}\phi)
\end{align}
This newly introduced model is referred to as the Cauchy-type damage model. 

\begin{figure}[!htb]
\centering
\includegraphics[width=8cm]{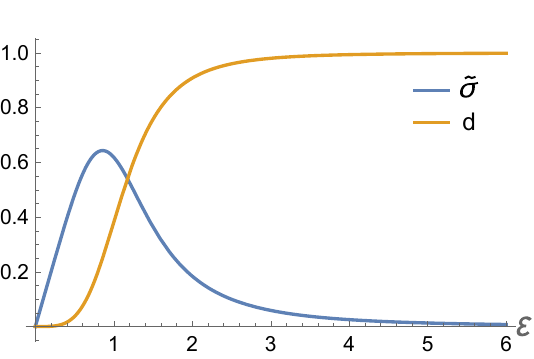}
\caption{Damage model based on Cauchy distribution in 1D with $G/\ell=1, k=1$: damage $d$ and effective stress $\tilde{\sigma}=\frac{ \sigma}{1+\phi^2\frac{\pi ^2 \ell ^2}{4 G^2}}$.}
\label{fig:DamageCrack2}
\end{figure}

When $\phi$ is small, the Taylor series of $\psi$ is
\begin{align}
\frac{2G}{\pi\ell} \arctan(\frac{\pi\ell}{2G}\phi)= \phi-\frac{\pi^2 \ell^2}{12 G^2} \phi^3+O(\phi^4),
\end{align}
which demonstrates that the model reduces to linear elasticity under undamaged conditions.
\begin{figure}[!htb]
\centering
\includegraphics[width=10cm]{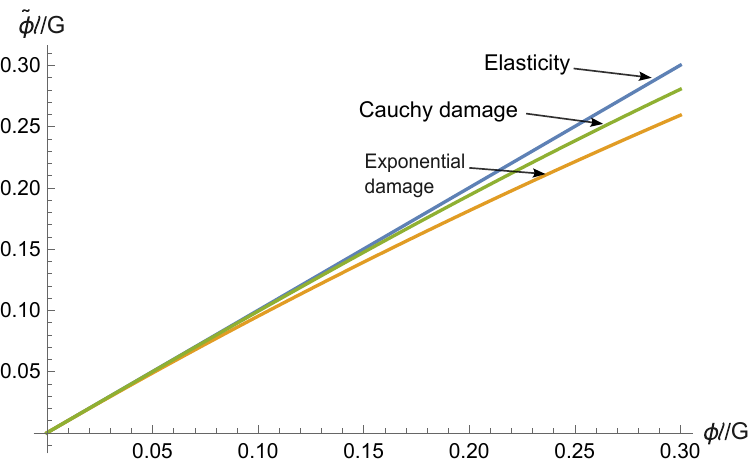}
\caption{Linear approximation of exponential damage model and Cauchy damage model when $\phi\ell/G<0.3$.}
\label{fig:LinearApproximationPhi}
\end{figure}
Figure \ref{fig:LinearApproximationPhi} shows the approximation of the Cauchy-type damage model on an elastic model. When the deformation energy satisfies $\phi\ell/G<0.3$, the maximal difference between damage model and elastic model is 6.55\%. In this sense, the damage model can account for the small deformation when damage is negligible.
In addition,
\begin{align}
\lim_{\phi\to +\infty}\frac{2G}{\pi\ell} \arctan(\frac{\pi\ell}{2G}\phi)=\frac{G}{\ell}
\end{align}
The variation of $\psi$ is
\begin{align}
\delta \psi=\frac{\delta \phi}{1+\phi^2\frac{ \pi ^2\ell ^2}{4 G^2}}
\end{align}
We employ the energy functional $E  =\int_\Omega \frac{2G}{\pi\ell} \arctan(\frac{\pi\ell}{2G}\phi^+)+\phi^--\bm b\cdot \bm u \, \ud V$. Through variational derivation of the energy functional, the governing equations are
\begin{align}
\nabla \cdot (\frac{ \bm \sigma^+}{1+\phi^{+2}\frac{\pi^2 \ell ^{2}}{4 G^2}}+\bm \sigma^-)+\bm b=\bm 0\label{eq:CauchyDamage}
\end{align}
For the convenience of post-processing, we introduce the damage variable defined by
\begin{align}
d=1-\frac{1}{1+\phi^{+2}\frac{\pi ^2 \ell ^2}{4 G^2}}
\end{align}
It can be easily verified that in 1D setting
\begin{align}
\int_{0}^{+\infty} \frac{ \sigma}{1+\phi^2\frac{\pi ^2 \ell ^2}{4 G^2}} \,\ud \varepsilon=\frac{G}{\ell}
\end{align}
Using relation $\frac{\partial \sigma}{\partial \varepsilon}=0$, the maximal stress is $\bm \sigma_{max}=\frac{3^{3/4}}{2} \sqrt{\frac{ G k}{\pi\ell}}$. The strain is $\varepsilon=\frac{2\cdot 3^{-1/4}}{\sqrt{\pi}}\sqrt{\frac{ G}{k \ell}}$. The corresponding damage is $d=1/4$.
In the setting of 1D with assumption $G/\ell=1, k=1$ and energy density $\phi=\frac 12 k\varepsilon^2$, we plot the graphs of damage function and effective stress in Figure \ref{fig:DamageCrack2}, which shows non-decreasing property of damage and the stress soften due to damage.
In order to determine the value of $\ell$, we employ the maximal stress determined from experiment. Similar to the routine in Eq.\ref{eq:Gellv}, the value of $G/\ell$ in current model can be solved from $\bm \sigma_{max}=\frac{3^{3/4}}{2} \sqrt{\frac{ G k}{\pi\ell}}$ as
\begin{align}
\frac{G}{\ell}=\frac{4 \pi}{3^{3/2}} \frac{\sigma_{max}^2}{k}. \label{eq:Gellv1}
\end{align}

\subsection{Damage model based on radical distribution}
Using Eq. \ref{eq:newProb1} and replacing the random variable $x$ with $\phi\ell/G$, the transformation is expressed as:
\begin{align}
\frac{x}{(1+x^{1/n})^n}\,\to \, \frac{G}{\ell} \frac{\phi\ell/G}{(1+(\phi\ell/G)^{1/n})^n}= \frac{\phi}{(1+(\phi\ell/G)^{1/n})^n}
\end{align}
Accordingly, the fracture energy density is defined as:
\begin{align}
\psi= \frac{\phi^+}{(1+(\phi^+\ell/G)^{1/n})^n}
\end{align}
and the corresponding energy functional is given by:  
\begin{align}
E   = \int_\Omega \left( \frac{\phi^+}{(1+(\phi^+\ell/G)^{1/n})^n} + \phi^- - \bm b  \cdot \bm{u} \right) \, \ud V.
\end{align}
Applying variational principles and integration by parts yields the governing equations for an elastic solid embedded with this damage model:  
\begin{align}
\nabla \cdot \Big(\frac{\bm \sigma^+}{(1+(\phi^+\ell /G )^{1/n})^{n+1}}+\bm \sigma^-\Big)+\bm b=\bm 0. \label{eq:gvm}
\end{align}
This damage model, represented by Eq.~\ref{eq:gvm}, aligns with the variational damage model described in Ref.~\cite{Ren2024vdm1}. The damage variable, used for post-processing purposes, is defined as:  
\begin{align}
d = 1 - \frac{1}{(1+(\phi^+\ell/G)^{1/n})^{n+1}}.
\end{align}

\begin{figure}[!htb]
\centering
\subfigure[]{
\includegraphics[width=0.45\textwidth]{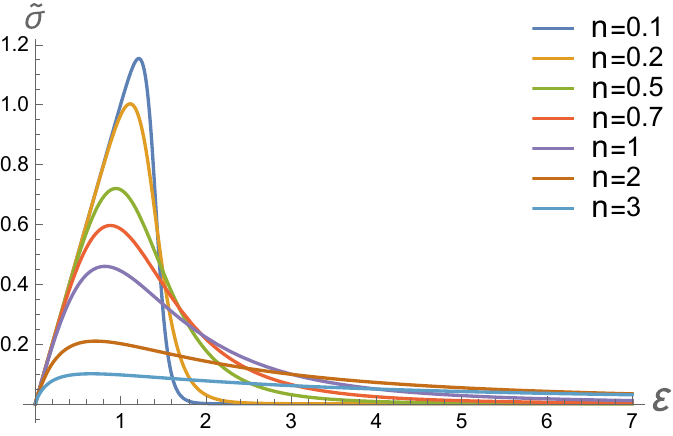}
\label{fig:sigmaEffectiveMs1}}
\subfigure[]{
\includegraphics[width=0.45\textwidth]{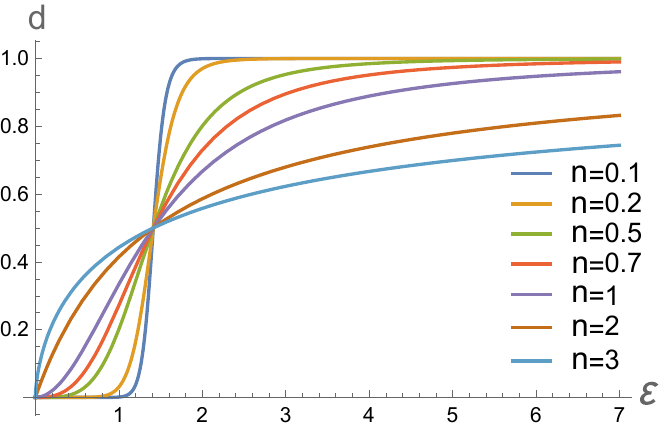}
\label{fig:damageEffectiveM0b}}
\caption{Damage model by Eq.\ref{eq:gvm} in 1D: (a) the effective stress versus strain; (b) the damage versus strain.}
\label{fig:EffectiveMs0}
\end{figure}

For the one-dimensional case, assuming $(G/\ell = 1, k = 1$), and an energy density $\phi = \frac{1}{2}k\varepsilon^2$, the evolution of effective stress and damage is illustrated in Figure \ref{fig:EffectiveMs0}.

In this formulation, the crack length scale $\ell$ remains undetermined. To estimate $\ell$, the maximum stress value $\bm \sigma_{max}$ obtained experimentally is used. Specifically, $(G/\ell)$ is derived from $\bm \sigma_{max}$ as:  
\begin{align}
\frac{ G}{\ell}= \frac{(1+\frac{n}{n+2})^{2(n+1)}}{2(\frac{n}{n+2})^{n}}\frac{\sigma_{max}^2}{k}.
\end{align}

\subsection{Counterexample based on Chi-square distribution}

In this section, we have derived four damage models from distinct probability distributions. More CDF-based damage models based on other feasible probability distributions are supplemented in \ref{sec:probDamage}.

Although the pool of available probability distributions is vast, not all can be effectively applied to damage modeling by simple algebraic replacement. A common feature among these applicable models is that the Taylor series expansion of the cumulative distribution function at $\phi=0$ includes linear terms of $\phi$. To illustrate this characteristic, we take the Chi-square distribution and replace the random variable with the strain energy density as follows: 
\begin{align}
\frac{1}{2 \Gamma\left(\frac{n}{2}\right)} \int_0^x\left(\frac{t}{2}\right)^{\frac{n}{2}-1} e^{-\frac{t}{2}} \ud t\to \frac{G}{\ell}\frac{1}{2 \Gamma\left(\frac{n}{2}\right)} \int_0^{2\phi\ell/G}\left(\frac{t}{2}\right)^{\frac{n}{2}-1} e^{-\frac{t}{2}} \ud t\notag
\end{align}

Based on this transformation, the fracture energy density is defined as:
\begin{align}
\psi=\frac{G}{\ell}\frac{1}{2 \Gamma\left(\frac{n}{2}\right)} \int_0^{2\phi\ell/G}\left(\frac{t}{2}\right)^{\frac{n}{2}-1} e^{-\frac{t}{2}} \ud t
\end{align}

The stress tensor influenced by the damage is computed using the relation:
\begin{align}
\tilde{\bm \sigma}=\frac{\partial \psi}{\partial\bm\varepsilon}=\frac{\partial \psi}{\partial\phi} \frac{\partial \phi}{\partial\bm\varepsilon}=\frac{1}{\Gamma\left(\frac{n}{2}\right)} \left(\frac{\phi\ell}{G}\right)^{\frac{n}{2}-1} e^{-\frac{\phi\ell}{G}} \bm \sigma \label{eq:Chisigma}
\end{align}
where $\bm \sigma=\frac{\partial \phi}{\partial\bm\varepsilon}$ represents the stress corresponding to the original energy density $\phi$. 

\begin{figure}[!htb]
\centering
\includegraphics[width=9cm]{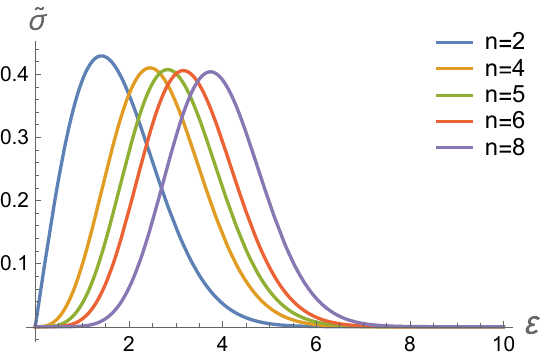}
\caption{Effective stress of the damage model based on Chi-square distribution in 1D.}
\label{fig:Chisigma}
\end{figure}

For a one-dimensional case under the assumptions $G/\ell=1$ and $k=1$, the strain energy density is given by $\phi=\frac{1}{2} k\varepsilon^2$. The evolution of stress, computed using Eq.~\eqref{eq:Chisigma}, is plotted in Figure~\ref{fig:Chisigma} for $n\in \{2,4,5,6,8\}$. The results indicate that the curve corresponding to $n=2$ successfully captures the stress-softening behavior observed in materials, while other values fail to adequately describe the linear elasticity phase in the regime of small deformations with negligible damage.

Interestingly, for $n=2$, the fracture energy density simplifies to:
\begin{align}
\psi=\frac{G}{\ell}\Big(1-\exp\big(-\frac{\ell}{G} \phi\big)\Big),\notag
\end{align}
which is equivalent to the exponential damage model. 

\section{Compactness and $\Gamma$-convergence of the CDF-based energy functional}
\newcommand{\eps}{\varepsilon}

\begin{theorem}[Compactness and $\Gamma$-convergence]\label{thm:Gamma}
Let  $\Omega\subset R^d$ (with $d=2,3$) be a bounded Lipschitz domain. Let $F:[0,\infty)\to[0,1]$
be a nondecreasing $C^1$ function such that $F(0)=0$, and there exists a critical value $s_c>0$ with
\begin{align}
F(s)=1\quad\text{for all }s\ge s_c,
\end{align} 
so that $F$ saturates rapidly near $s_c$. For each $\lambda>0$, set
\begin{align}
\psi_\lambda(t) := \int_0^t [1-F_\lambda(s)] \ud s\,=\frac{1}{\lambda}\int_0^{\lambda t} [1-F(s)] ds,\qquad t\ge 0.
\end{align} 
with $F_\lambda(s) := F(\lambda s)$.
The scaled CDF $F_\lambda(s)$ satisfies the saturation property 
\begin{align}
1-F_\lambda(s)=0 \quad \mbox{for  } s\ge s_c/\lambda.
\end{align} 

Consider the family of energy functionals $E_\lambda:\mbox{SBV}(\Omega)\to[0,\infty]$ given by
\begin{align}
E_\lambda(u) \;=\; \int_\Omega \psi_\lambda \Bigl(\frac{|\nabla^a u|^2}{2}\Bigr)\,\ud V, 
\end{align}
where $\nabla^a u$ denotes the absolutely continuous part of $Du$ (and $E_\lambda(u)=+\infty$ if $u$ has any Cantor part). Assume moreover that the functions $u$ are uniformly bounded in $L^\infty(\Omega)$, or have prescribed bounded boundary data. Then the following statements hold as $\lambda\to+\infty$:

\textbf{Compactness}:   If a sequence $(u_\lambda)\subset \mbox{SBV}(\Omega)$ satisfies
\begin{align}
\sup_\lambda\Bigl\{\|u_\lambda\|_{L^\infty(\Omega)} + E_\lambda(u_\lambda)\Bigr\}<+\infty, 
\end{align}
then, up to a subsequence, $u_\lambda\to u$ strongly in $L^1(\Omega)$ for some $u\in \mbox{SBV}(\Omega)$ whose approximate gradient vanishes almost everywhere; that is, $\nabla^a u=0$. In particular, the limit $u$ is (up to null sets) a piecewise constant function.

\textbf{$\Gamma$-convergence}:  The functionals $E_\lambda$ $\Gamma$-converge with respect to the strong $L^1(\Omega)$-topology to the limiting functional
\begin{align}
E_\infty(u) \;=\; \begin{cases} 0,& u\in \mbox{SBV}(\Omega), \nabla^a u=0\text{ a.e. in }\Omega,\\ +\infty,& \text{otherwise}. \end{cases} 
\end{align}
In other words, any $L^1$-limit of finite-energy sequences must have zero approximate gradient, and conversely any piecewise-constant (i.e. $\nabla^a u=0$) function can be approximated in $L^1$ by a sequence whose energies $E_\lambda$ vanish as $\lambda\to\infty$.
\end{theorem}

\noindent\textbf{Proof.}

\textbf{Compactness:} Since $F(s)=1$ for $s\ge s_c$, each integrand $\psi_\lambda(t)=\int_0^t [1-F(\lambda t)]dt$ saturates at the value 1 when $t\ge s_c/\lambda$. 
Hence
\begin{align}
0\le \psi_\lambda\Bigl(\tfrac{|\nabla^a u_\lambda|^2}{2}\Bigr)\le t,
\end{align} 
and $E_\lambda(u_\lambda)\le t|\Omega|$ uniformly. More importantly, a uniform bound $E_\lambda(u_\lambda)\le C$ forces the measure of the set
to be $\le C$. As $\lambda\to\infty$, the threshold $s_c/\lambda\to0$, so in the limit almost every point must satisfy $|\nabla^a u(x)|=0$. A more detailed argument using De Giorgi–Ambrosio compactness for $SBV$ (see e.g.\ \cite{ambrosio2000functions,dal1993introduction}) shows that $(u_\lambda)$ is relatively compact in $L^1$, and any limit $u$ has no absolutely continuous gradient part, i.e.\ $\nabla^a u=0$ a.e.

\textbf{$\Gamma$-liminf:} Take any sequence $u_\lambda\to u$ in $L^1(\Omega)$. If $\nabla^a u\neq0$ on a set of positive measure, then for large $\lambda$ the gradient $|\nabla^a u_\lambda|$ must be nonzero on a positive-measure set as well (by lower-semicontinuity of the gradient part), so $\psi_\lambda(|\nabla^a u_\lambda|^2/2)$ saturates to $t$ there. 

Hence $\liminf_{\lambda}E_\lambda(u_\lambda)\ge||\nabla^a u||>0$, forcing $E_\infty(u)=+\infty$ by definition. Conversely, if $\nabla^a u=0$ a.e., then $u$ is piecewise constant up to a jump set. One constructs a recovery sequence by ``smoothing'' each jump of $u$ over a thin layer of width $\eps (\lambda) \to 0$, with gradient of size $O(1/\eps)$ so that $|\nabla u_\lambda|^2\sim O(1/\eps^2)$ and $\int_{I_\eps}  \psi_\lambda\Bigl(\tfrac{|\nabla^a u_\lambda|^2}{2}\Bigr) \approx 0\cdot $(width) in this layer (the integral of $1-F_\lambda$ is zero once the gradient exceeds the threshold; hence $\psi_\lambda$ saturates at the finite value $G:=\int_{0}^{\infty}[1-F(s)]ds$. The recovery sequence is built so that each layer costs $\approx G\,\eps \to 0$). The contribution of each jump of (hyper)surface area $A$ to $E_\lambda$ is then $1\times(A,\eps)\approx (A,\eps)$, which tends to $0$ as $\eps\to 0$. Outside these layers, $u_\lambda$ is constant and contributes no energy. Hence $E_\lambda(u_\lambda)\to 0=E_\infty(u)$, proving the $\Gamma$-limsup inequality.
This establishes the claimed $\Gamma$-convergence.

\section{Existence of Rate-Independent Evolution of the CDF-based energy functional}

\begin{theorem}[Quasi-static evolution under loading]\label{thm:RI-evol}
In the setting of Theorem \ref{thm:Gamma}, fix a parameter $\lambda>0$ and consider the (time-dependent) energy functional
\begin{align}
\mathcal{E}(t,u) \;=\;\int_\Omega \int_0^{\tfrac{|\nabla^a u|^2}{2}}[1- F (\lambda s)]\ud s\,\ud V, 
\end{align}
defined for $u\in SBV(\Omega)$ that satisfy a given time-dependent Dirichlet boundary condition $u|_{\partial\Omega}=g(t)$ with $g\in C^1([0,T];H^{1/2}(\partial\Omega))$. 

Define a dissipation distance $\mathcal{D}:SBV(\Omega)\times SBV(\Omega)\to[0,\infty]$ by
\begin{align}
\mathcal{D}(u_0,u_1) = \begin{cases}
 \int_{J_{u_1}\setminus J_{u_0}}\gamma\bigl(|[u_1]|\bigr)\,d H^{d-1}, &\text{if }J_{u_0}\subseteq J_{u_1},\\[1.2em] +\infty,&\text{otherwise}, 
 \end{cases} 
\end{align}
where $J_u$ is the jump set of $u$, $[u_1]$ the jump height, and $\gamma:[0,\infty)\to(0,\infty)$ is a continuous cost-density with $\gamma(0)=0$ (e.g. a positive constant). This $\mathcal{D}$ enforces irreversibility (no healing of cracks) and satisfies the triangle inequality.

Assume an initial condition $u(0)\in \mbox{SBV}(\Omega)$ is given and is stable at $t=0$, i.e. $\mathcal{E}(0,u(0)) \;\le\; \mathcal{E}(0,v) + \mathcal{D}(u(0),v) \quad\text{for all admissible }v\text{ with }v|_{\partial\Omega}=g(0)$. 

Then there exists an energetic solution $u:[0,T]\to \mbox{SBV}(\Omega)$ of the rate-independent system. That is, $u(t)$ satisfies for all $t\in[0,T]$:
\begin{itemize}
\item[(S)] (Global stability.)
\begin{align}
\mathcal{E}(t,u(t)) \;\le\; \mathcal{E}(t,v) + \mathcal{D}\bigl(u(t),v\bigr) \quad\text{for every admissible }v,\;v|_{\partial\Omega}=g(t). 
\end{align}
\item[(E)] (Energy balance.)
\begin{align}
\mathcal{E}(t,u(t)) + \mathrm{Diss}_{\mathcal{D}}\bigl(u;[0,t]\bigr) \;=\; \mathcal{E}(0,u(0)) + \int_0^t \partial_s \mathcal{E}(s,u(s))\,ds, 
\end{align}
where $\mathrm{Diss}{\mathcal{D}}(u;[0,t])$ is the total $\mathcal{D}$-dissipation accumulated by the path $u(\cdot)$ over $[0,t]$.
\end{itemize}
In particular, the jump set $J{u(t)}$ grows monotonically (no healing), and the above equalities hold for almost every $t$.
\end{theorem}

\noindent\textbf{Proof.}

We outline the standard time-discretization argument (cf.\ \cite{mielke1999mathematical,mainik2005existence}). By Theorem \ref{thm:Gamma}, each sublevel set of $\mathcal{E}(t,\cdot)$ is sequentially $L^1$-compact in SBV (since high gradients are saturated and promote compactness via piecewise-constancy). Moreover $\mathcal{E}(t,u)$ is lower semicontinuous in $L^1$ and depends continuously on $t$ through $g(t)$, so $\partial_t\mathcal{E}(t,u)$ is bounded. The dissipation $\mathcal{D}(u_0,u_1)$ is lower semicontinuous in $u_1$ (under $L^1$-convergence) and satisfies the triangle inequality (by additivity of new crack measure). Hence all abstract hypotheses for energetic solutions are met (cf. \cite{mielke1999mathematical,mainik2005existence}).

Concretely, for a partition $0=t_0<t_1<\cdots<t_N=T$ one constructs inductively minimizers $u^k\in \mbox{SBV}(\Omega)$ of the incremental problem
\begin{align}
u^k \;\in\; \arg\min_{\,v}\Bigl\{\;\mathcal{E}(t_k,v)\;+\;\mathcal{D}(u^{k-1},v)\Bigr\}, 
\end{align}
with $u^0$ given. Existence of each minimizer follows by the direct method: $\mathcal{E}$ is $L^1$-lower semicontinuous with compact sublevels (by the compactness above) and $\mathcal{D}(u^{k-1},v)$ imposes a constraint $J_{u^{k-1}}\subset J_v$ that is preserved under $L^1$-limits, so admissible competitors form a closed set. Hence a minimizer exists for each step.

Passing to the time-continuous limit as the mesh of the partition goes to zero yields, via a Helly-type compactness argument, a limit trajectory $u(t)$ which satisfies the discrete stability at all rational times and the discrete energy inequality. By standard arguments (cf.\cite{mainik2005existence}), this limit in fact satisfies the global stability (S) and the energy balance (E) in the statement, with total dissipation equal to the sum of increments $\sum_k\mathcal{D}(u^{k-1},u^k)$. In particular one uses the compatibility condition that $J_{u^k}\subseteq J_{u^{k+1}}$, so there is no healing; the irreversibility is enforced by $\mathcal{D}(u_{k-1},u_k)=+\infty$ if $J_{u_k}$ fails to contain $J_{u_{k-1}}$.

Thus we obtain an energetic solution $u(t)$ with the desired properties (see also \cite{mielke2003rate,maso2005quasistatic} for similar constructions) by existence of minimizers at each step and compactness (incremental problem), and then by a Helly selection to pass to a continuous-time solution satisfying (S) and (E).

\textbf{Remark:} The assumptions on $F$ (rapid saturation) and on $\gamma$ ensure that $E_\lambda$ is coercive in the SBV-sense (i.e.\ pushes the absolutely continuous gradient to zero) while $\mathcal{D}$ penalizes the creation of new jump surface. One recovers a purely rate-independent fracture behavior: the solution $u(t)$ remains piecewise constant except where forced to evolve by the loading $g(t)$, and any jump created carries a finite cost in the energy balance.

\section{Illustrative Example: Single-Edge-Notched Plate under Tension}
\begin{figure}[htbp]
\centering
\subfigure[]{\includegraphics[width=0.3\textwidth]{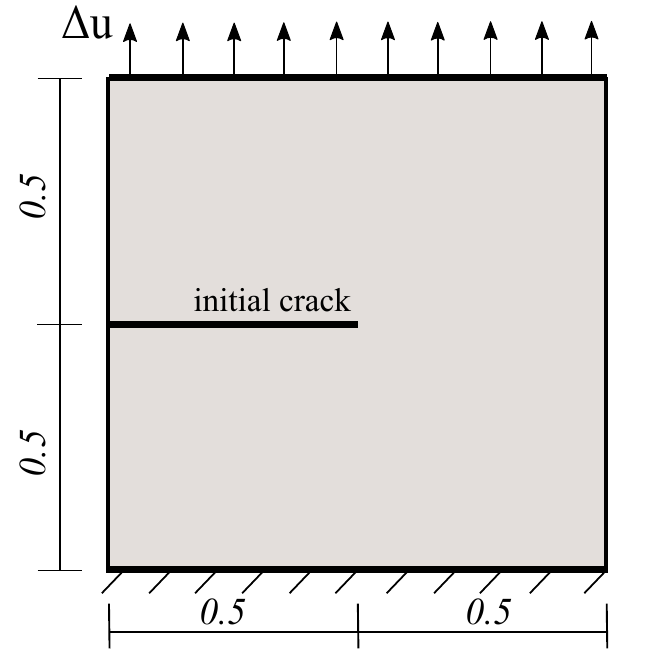}}
\subfigure[]{\includegraphics[width=0.27\textwidth]{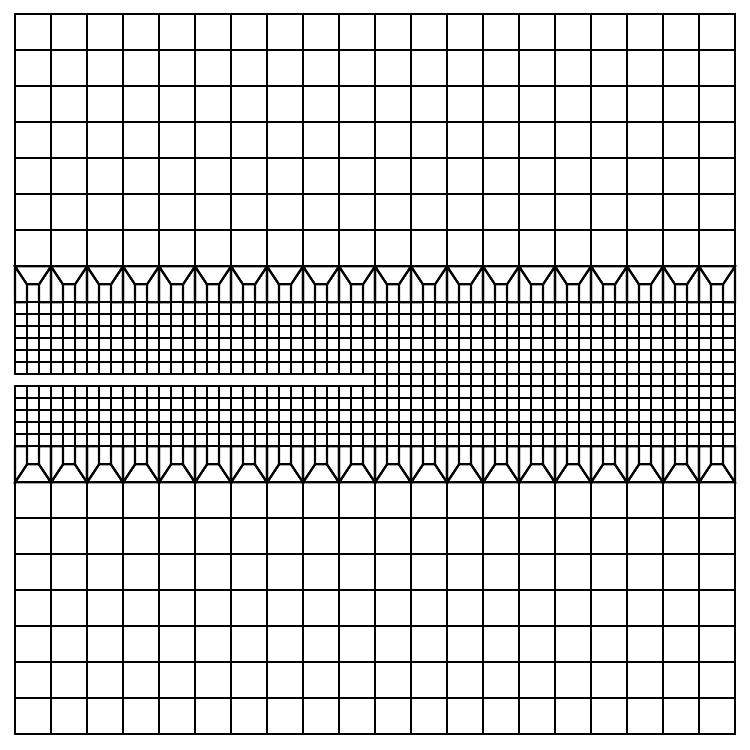}}
\subfigure[]{\includegraphics[width=0.27\textwidth]{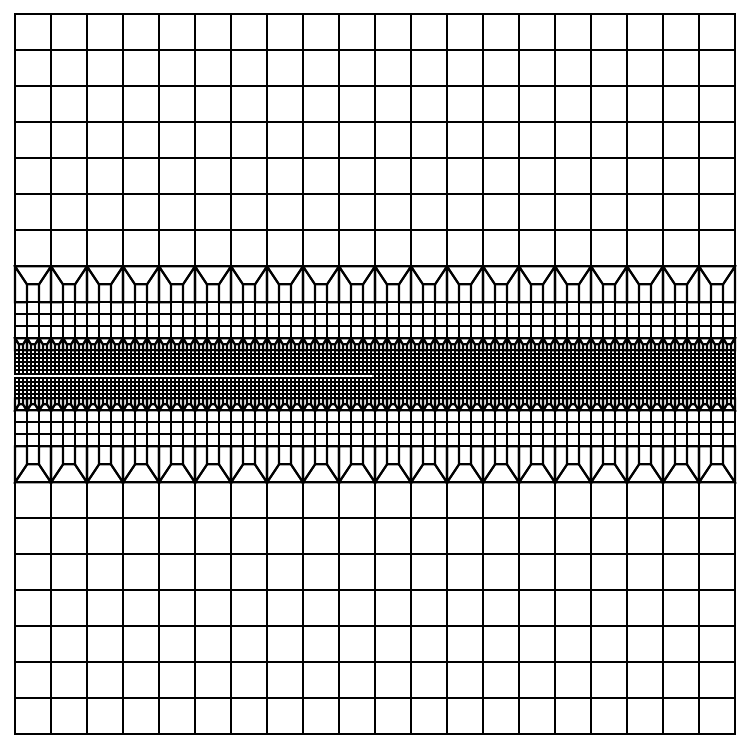}}
\caption{Single-edge-notched plate: (a) geometry and boundary conditions (unit: mm) ; (b) FEM coarse mesh; (c) FEM refined mesh.}
\label{fig:PlateSetUp}
\end{figure}

The following example summarizes one benchmark also studied in \cite{RenCDFfem2025} but now used to illustrate the mathematical formulation. We consider a standard single-edge-notched tensile (SENT) specimen under plane strain conditions  in two dimensions. The computational domain, boundary conditions, and geometric setup are shown in Figure~\ref{fig:PlateSetUp}(a). The bottom edge is fixed in both horizontal and vertical directions, while a uniform vertical displacement is applied along the top edge. No body forces are applied. The material parameters are selected to represent a typical metal: Young’s modulus \(E = 210\;\mathrm{kN/mm}^2\), Poisson’s ratio \(\nu = 0.3\), and critical energy release rate \(G_c = 2.7 \times 10^{-3}\;\mathrm{kN/mm}\). The fracture energy density is computed directly at integration points.

Figure~\ref{fig:PlateSetUp}(b,c) display the discretisation used in the simulation. Two quadrilateral (Q4) mesh resolutions are considered: a coarse mesh with 3520 elements and a refined mesh with 9280 elements. Both meshes are locally refined in the vicinity of the notch to better capture crack initiation and propagation. The corresponding effective element sizes are \(\Delta x = 1/120\;\mathrm{mm}\) and \(\Delta x = 1/360\;\mathrm{mm}\), respectively.

Five representative CDF-based damage models are tested: exponential, Cauchy-type, logistic-type, half-normal-type, and Gudermannian-type. Each model defines a distinct softening response and therefore requires a distinct internal length scale \(\ell\). For numerical consistency, the characteristic lengths are set to  $\ell = (1.0462,\, 1.0955,\, 1.4203,\, 1.5071,\, 1.3505)\cdot \Delta x,$ respectively, for the five models listed above.

\begin{figure}[htbp]
\centering
\subfigure[Coarse mesh]{\includegraphics[width=0.42\textwidth]{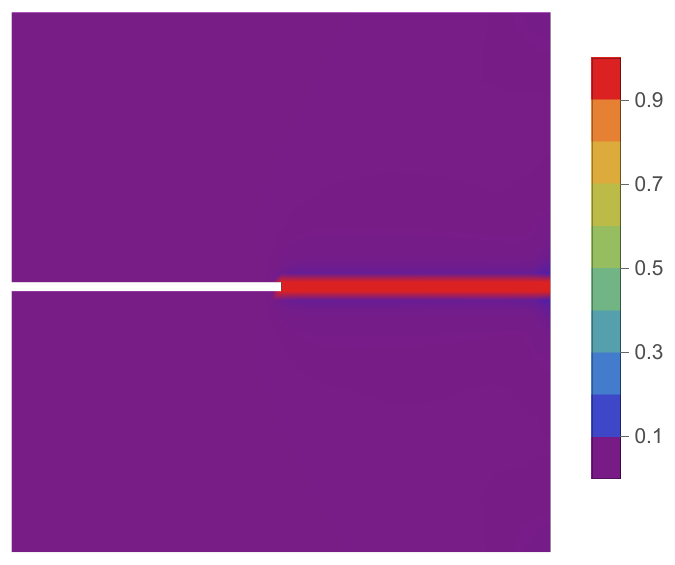}}
\quad
\subfigure[Refined mesh]{\includegraphics[width=0.42\textwidth]{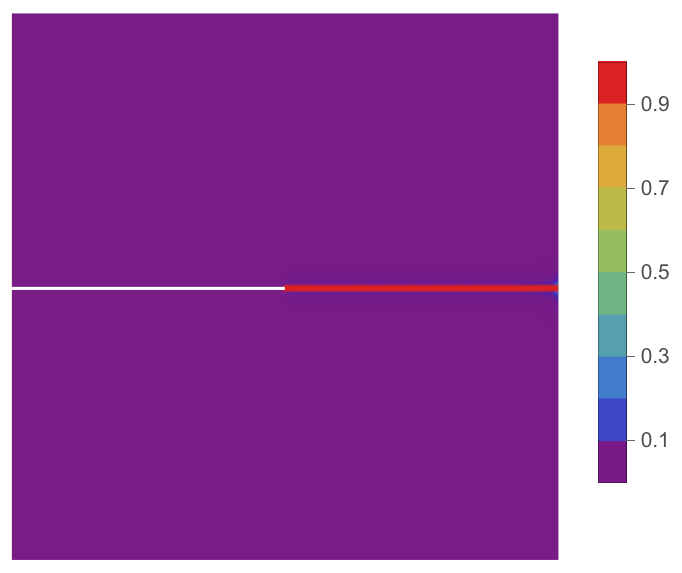}}
\caption{Damage distribution for the logistic-type model: (a) coarse mesh; (b) refined mesh.}
\label{fig:Damage2}
\end{figure}

\begin{figure}[htbp]
\centering
\subfigure[Coarse mesh]{\includegraphics[width=0.42\textwidth]{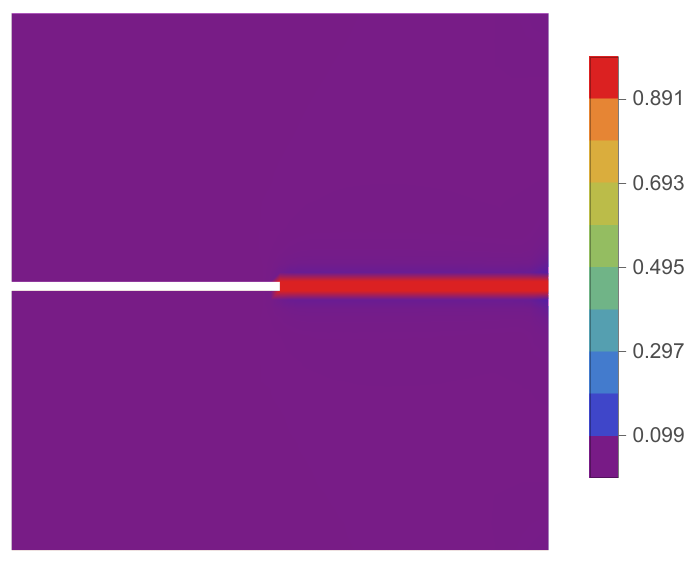}}
\quad
\subfigure[Refined mesh]{\includegraphics[width=0.42\textwidth]{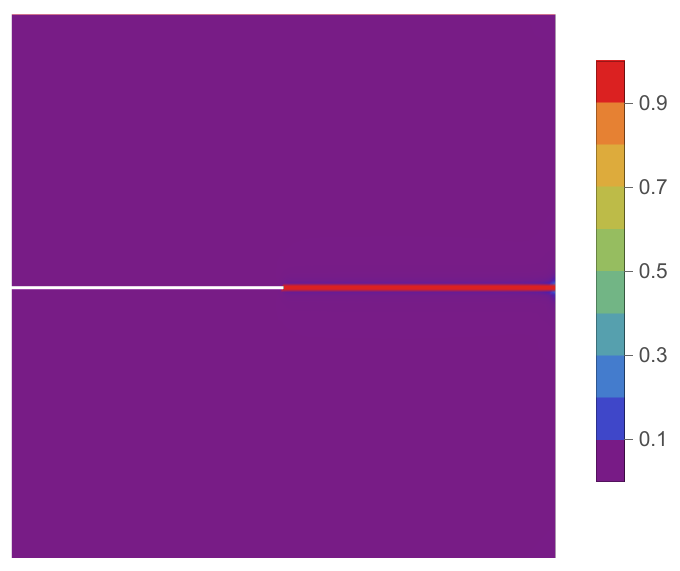}}
\caption{Damage distribution for the half-normal-type model: (a) coarse mesh; (b) refined mesh.}
\label{fig:Damage5}
\end{figure}

Figures~\ref{fig:Damage2} and \ref{fig:Damage5} display the predicted damage fields for the logistic-type and half-normal-type models, respectively, on both mesh resolutions. In both cases, the models produce a sharply localized crack pattern aligned with the expected path of tensile failure. The crack trajectory remains mesh-insensitive, confirming the robustness of the formulation.

\begin{figure}[htbp]
\centering
\includegraphics[width=0.55\textwidth]{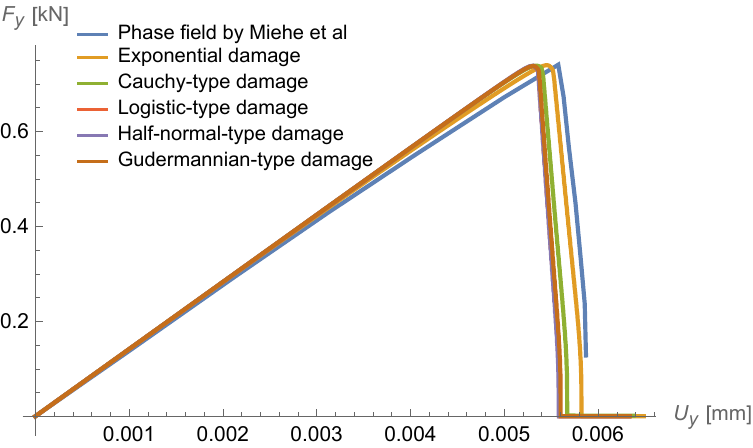}
\caption{Load–displacement curves of the SENT plate for five CDF-based damage models (coarse mesh). Comparison with the phase-field model of Miehe \etal~\cite{miehe2010thermodynamically}.}
\label{fig:Damage1-9LCcoarse}
\end{figure}

The load–displacement responses for all five damage models using the coarse mesh are plotted in Figure~\ref{fig:Damage1-9LCcoarse}. All curves show good agreement with the phase-field fracture results of Miehe \etal \cite{miehe2010thermodynamically}, confirming the physical fidelity of the CDF-based models.

\begin{figure}[htbp]
\centering
\subfigure[Exponential]{\includegraphics[width=0.45\textwidth]{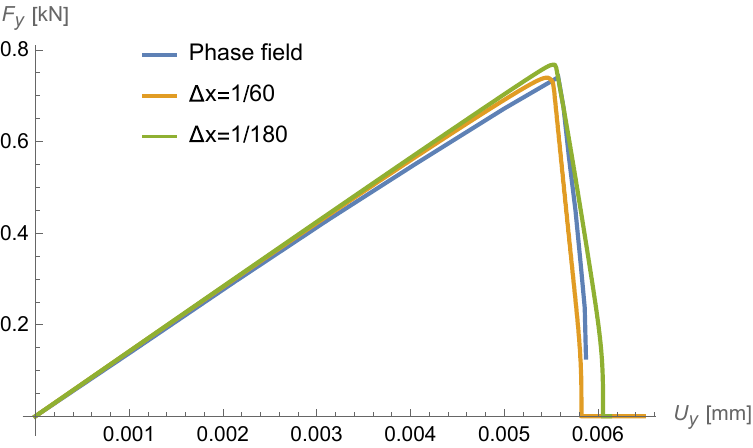}}
\subfigure[Cauchy-type]{\includegraphics[width=0.45\textwidth]{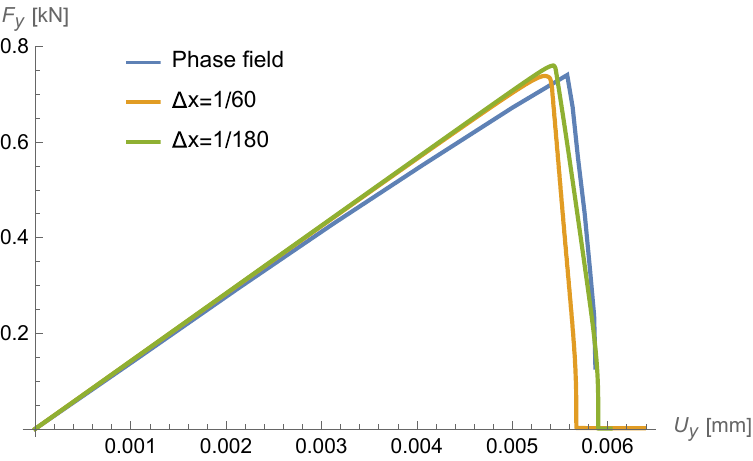}}\\
\subfigure[Logistic-type]{\includegraphics[width=0.45\textwidth]{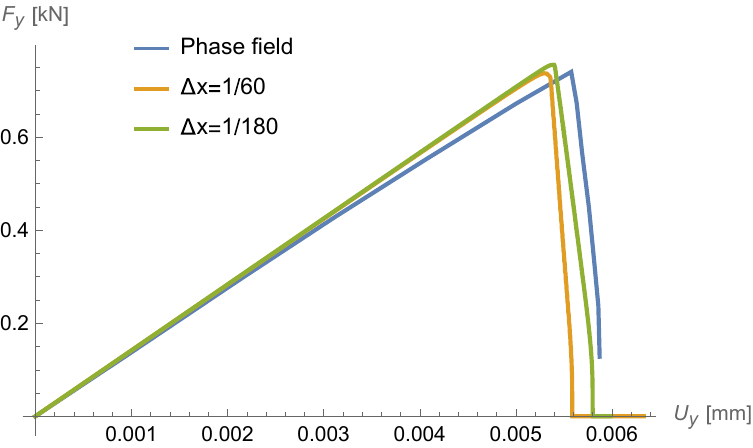}}
\subfigure[Half-normal-type]{\includegraphics[width=0.45\textwidth]{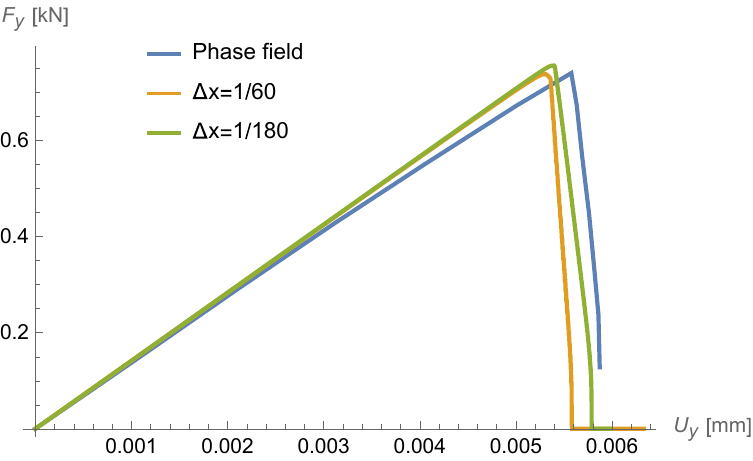}}
\caption{Mesh sensitivity study: load–displacement curves for four CDF-based models on coarse vs. refined meshes. Comparison with phase-field benchmark results from \cite{miehe2010thermodynamically}.}
\label{fig:PFvsVD}
\end{figure}

Finally, the influence of mesh resolution on the peak load is examined for four selected models. Figure~\ref{fig:PFvsVD} shows that both coarse and refined meshes yield nearly identical peak loads, indicating good convergence. A full 3-D validation and mesh convergence study are reported in \cite{RenCDFfem2025}. 

\section{Conclusion remarks}

We have provided a variational foundation for various analytic softening laws obtained by integrating the complements of well-known probability distributions. The core findings are:

1.	\textbf{Admissibility Theorem} Any CDF with $F(0)=0$ and $F(\infty)=1$ yields a damage map that is monotone, bounded, and dissipative, ensuring thermodynamic consistency for the corresponding stored-energy density.

2.	\textbf{$\Gamma$-convergence to Griffith fracture} For every sub-family with finite saturation energy $G$ the  energy sequence $\Gamma$-converges to the Griffith functional.

3. Explicit recovery sequences and a compactness lemma establish both the liminf and limsup inequalities.

4.	\textbf{Well-posedness} The CDF-based damage functional admits an energetic solution that evolves quasi-statically under a rate-independent loading path.

Collectively, these results supply a bridge between probability-inspired softening laws and classical fracture mechanics, clarifying when and how diffuse damage degenerates into a sharp-crack model. Although the present work is focused on theoretical analysis, the proposed model has been implemented in a finite element framework and applied to a range of engineering fracture problems. The companion paper \cite{RenCDFfem2025} presents detailed numerical studies, including convergence tests and comparisons to established methods, thereby validating the applicability of the theory. Future work will generalise the $\Gamma$-convergence proof to tensorial damage and ductile damage.

\section*{Acknowledgments}
The research is supported by the Fundamental Research Funds for the Central Universities.

\section*{Author Contributions Statement}
\noindent
\textbf{H. Ren}: Conceptualization, Methodology, Formal analysis, Mathematical proofs, Writing – original draft, Writing – review \& editing, Visualization.

\vspace{1em}
\section*{Ethics Approval and Consent to Participate}
\noindent
Not applicable. This study does not involve human participants, animals, or identifiable personal data.

\appendix

\section{Several conventional distributions}\label{sec:probdata}
\subsection{Cauchy distribution}
The Cauchy distribution is a heavy-tailed distribution with undefined mean and variance, characterized by its lack of moments and a probability density function with a peak at its location parameter but with very slowly decaying tails. It is often used to illustrate phenomena where extreme values are more likely than in normal distributions.
The probability density function of the Cauchy distribution, centered at zero, is given by:
\begin{align}
f(x; 0, \gamma) = \frac{1}{\pi} \cdot \frac{\gamma}{x^2 + \gamma^2},
\end{align}
where $\gamma > 0$ is the scale parameter, determining the width of the peak.
The cumulative distribution function of the Cauchy distribution is expressed as:
\begin{align}
F(x; 0, \gamma) = \frac{1}{\pi} \arctan\left(\frac{x}{\gamma}\right) + \frac{1}{2}. \label{eq:Caudis}
\end{align}
A comprehensive explanation of the Cauchy distribution can be found in Ref \cite{walck1996hand}.

\subsection{Logistic distribution}
The logistic distribution resembles the normal distribution but has heavier tails, characterized by a symmetric, S-shaped CDF. It is often used in logistic regression and modeling growth processes, where its heavier tails provide robustness against outliers compared to the normal distribution.
The cumulative distribution function of the logistic distribution is defined as:
\begin{align}
F(x) = 1 - \frac{1}{1 + e^x}, \quad x \in (-\infty, +\infty). \label{eq:Logdis}
\end{align}
The probability density function is obtained by differentiating the CDF, and is given by:
\begin{align}
f(x) = F'(x) = \frac{2 e^x}{(1 + e^x)^2}, \quad x \in (-\infty, +\infty).
\end{align}
For additional insights into the logistic distribution, please consult Ref \cite{walck1996hand}.

\subsection{Half-normal distribution}
The half-normal distribution is a special case of the normal distribution restricted to non-negative values, resulting from taking the absolute value of a normally distributed variable with mean zero and a given standard deviation. It is used to model absolute deviations and is characterized by its single peak at zero and a rapidly decreasing tail \cite{tsagris2014folded}.
In this distribution, the CDF is expressed using the error function, $\erf(x)$, which is also known as the Gauss error function. The error function is defined as:
\begin{align}
\erf(x) = \frac{2}{\sqrt{\pi}} \int_0^x e^{-t^2} \, dt.
\end{align}
For a random variable $x$ following the half-normal distribution, the CDF is given by:
\begin{align}
F(x) = \erf(x) = \frac{2}{\sqrt{\pi}} \int_0^x e^{-t^2} \, dt, \quad x \in [0, +\infty). \label{eq:errordis}
\end{align}
The corresponding probability density function is obtained by differentiating the CDF:
\begin{align}
f(x) = F'(x) = \frac{2 e^{-x^2}}{\sqrt{\pi}}.
\end{align}
One of the notable properties of the half-normal distribution is its moments. The $m$-th raw moment of the distribution is given by:
\begin{align}
\mathbb E(x^m) = \int_0^{+\infty} f(x) x^m \, \ud x  = \frac{\Gamma \left( \frac{m+1}{2} \right)}{\sqrt{\pi}}, \quad \text{if } \Re(m) > -1.
\end{align}
Here, $\Gamma(\cdot)$ is the Gamma function, which generalizes the factorial function to non-integer values. The reader is encouraged to explore Ref \cite{walck1996hand} for a thorough discussion of the half-normal distribution.

\subsection{Chi-square distribution}
The Chi-square distribution is a distribution of the sum of the squares of $n$ independent standard normal random variables, often used in hypothesis testing, particularly in tests of independence and goodness of fit.
Mathematically, the Chi-square distribution is given by
\begin{align}
f(x ; n)=\frac{\left(\frac{x}{2}\right)^{\frac{n}{2}-1} e^{-\frac{x}{2}}}{2 \Gamma\left(\frac{n}{2}\right)}
\end{align}
where the variable $x \geq 0$ and the parameter $n$, the number of degrees of freedom, is a positive integer.
The cumulative, or distribution, function for a Chi-square distribution with $n$ degrees of freedom is given by
\begin{align}
F(x) & =\frac{1}{2 \Gamma\left(\frac{n}{2}\right)} \int_0^x\left(\frac{t}{2}\right)^{\frac{n}{2}-1} e^{-\frac{t}{2}} d t=\frac{\gamma\left(\frac{n}{2}, \frac{x}{2}\right)}{\Gamma\left(\frac{n}{2}\right)}
\end{align}
where $\gamma(a,x)=\int_{0}^x t^{a-1} e^{-t} dt$ and $\Gamma(a,x)=\int_{x}^\infty t^{a-1} e^{-t} dt$ are the incomplete Gamma function with $\Re(a)>0$. Detailed information about the Chi-square distribution is available in Ref \cite{walck1996hand}.

\section{More new probability distributions}\label{sec:probdatanew}

\subsection{Rational distribution}
\begin{figure}[htp]
\centering
\subfigure[]{
\includegraphics[width=0.46\textwidth]{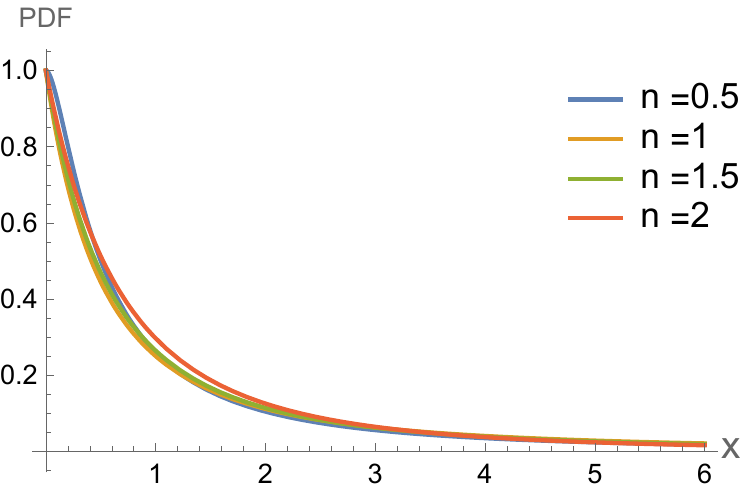}
\label{fig:raydisa2}}
\subfigure[]{
\includegraphics[width=0.46\textwidth]{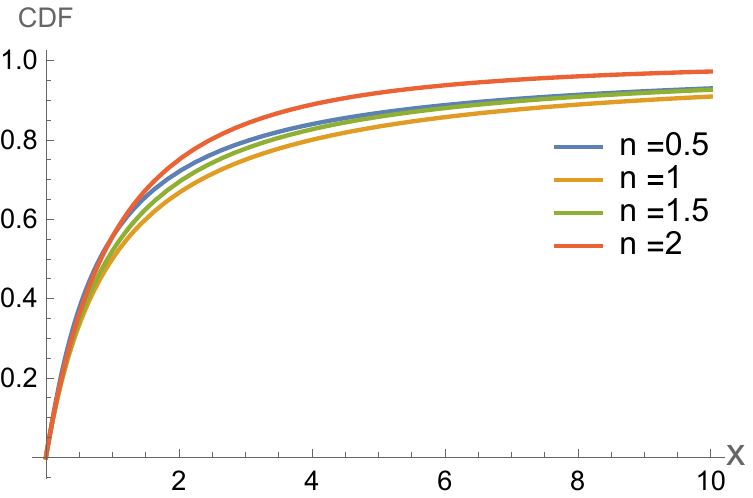}
\label{fig:raydisb2}}
\caption{Rational distribution: (a) probability density function (b) cumulative distribution function}
\label{fig:Piecewise3}
\end{figure}

By utilizing a rational polynomial function, we propose a new accumulated distribution function defined as:
\begin{align}
F(x;n)=\frac{x^2+x n^2}{(x+n)^2} ,\,x\in [0,+\infty),\label{eq:newProb3}
\end{align}
where $0.5\leq n\leq 2$ is the scale parameter of the distribution. It can be easily verified that $F(0;n)=0$, $\lim_{x\to +\infty} F(x;n)=1$ and $F(x;n)$ is a non-decreasing continuous function. We refer to this distribution as the rational distribution.
The probability density function corresponding to this CDF is:
\begin{align}
f(x;n)=F'(x;n)=\frac{n \left(n^2-n x+2 x\right)}{(n+x)^3}.\label{eq:pdf3}
\end{align}
Both the PDF and CDF of this rational distribution are illustrated in Figure \ref{fig:Piecewise3}. Notably, the distribution demonstrates relative insensitivity to the parameter $n$.
The $m$-th raw moment of the distribution is given by:
\begin{align}
\mathbb E(x^m)=\int_{0}^{+\infty} f(x;n) x^m \ud x =\pi m n^m (1+m-mn) \csc (\pi m)\mbox{ if }m<1.
\end{align}

\subsection{Gudermannian distribution}
The Gudermannian function, abbreviated $gd(x)$, is an odd function that is used in the inverse equations for the Mercator projection \cite{beyer1978crc}. Function $gd(x)$ in range $[0,+\infty)$ meets all the requirements of the CDF. Therefore, it is natural to derive Gudermannian distribution from the Gudermannian function. Let the CDF use the $gd(x)$ for random variable $x$:
\begin{align}
F(x)=gd(x)=\int_{0}^x \sech(\frac{\pi}{2}t)\, dt=\frac{4}{\pi} \arctan(\tanh(\frac{\pi}{4}x)), \, x\in [0,+\infty) \label{eq:gddis}
\end{align}
It is worthy noting that the distribution by Eq.\ref{eq:gddis} is different from the generalized Gudermannian distribution in Ref \cite{altun2019generalized}, which used Gudermannian function with the form of $\frac{2}{\pi}\arctan(\exp(x))$.
A compact PDF based on Eq.\ref{eq:gddis} can be derived as
\begin{align}
f(x)=\sech(\frac{\pi}{2} x),\, x\in [0,+\infty)
\end{align}
\begin{figure}[!htb]
\centering
\includegraphics[width=7cm]{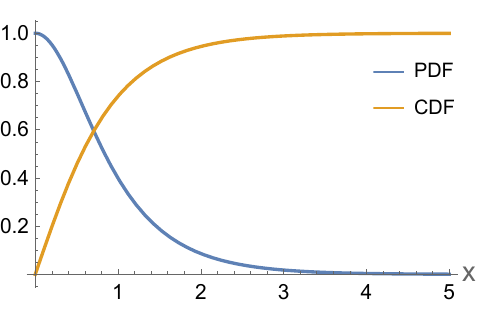}
\caption{Gudermannian distribution:PDF and CDF.}
\label{fig:guddis}
\end{figure}
The graphs of PDF and CDF of this distribution are plotted in Figure \ref{fig:guddis}. For this probability distribution, the $m$-th plain moment is
\begin{align}
\mathbb E(x^m)=\int_0^{+\infty} f(x) x^m=\frac{\Gamma (m+1)}{2^{2 m+1}} \big(\zeta (m+1,\frac{1}{4})-\zeta (m+1,\frac{3}{4})\big) \text{ if }\Re(m)>-1
\end{align}
where $\zeta(a,b)$ is the Zeta function.

\subsection{Continuous hypergeometric distribution}
\begin{figure}[htp]
\centering
\subfigure[]{
\includegraphics[width=0.46\textwidth]{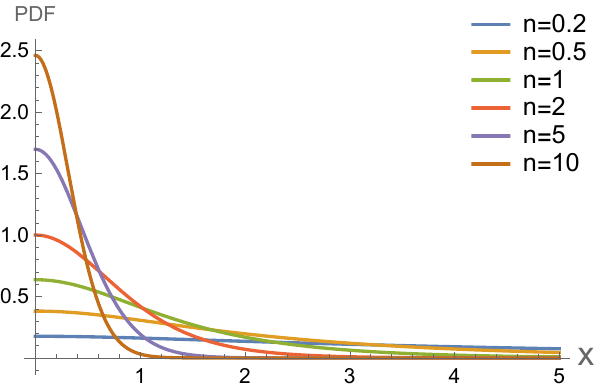}
\label{fig:raydisa}}
\subfigure[]{
\includegraphics[width=0.46\textwidth]{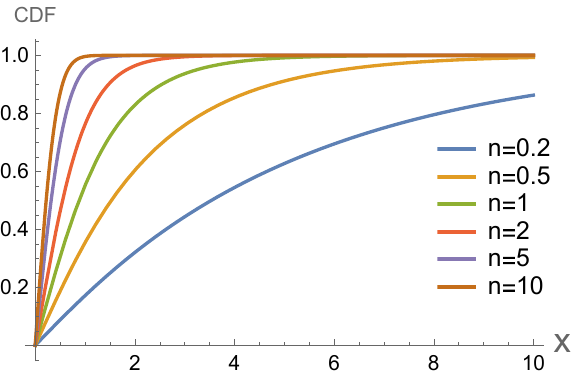}
\label{fig:raydisb}}
\caption{Continuous hypergeometric distribution: (a) probability density function (b) cumulative distribution function}
\label{fig:hfp}
\end{figure}

Hypergeometric distribution is usually referred to a discrete probability distribution that describes the probability of drawing a specific number of successes from a finite population without replacement, where the population contains a fixed number of successes and failures \cite{sibuya1979generalized,fog2008calculation}. In contrast to the discrete version, we propose a continuous hypergeometric distribution in this section.
The hypergeometric function is defined for $\|z\|\leq 1$ by the power series
\begin{align}
{}_{2}F_{1}(a,b;c;z)=\sum _{k=0}^{\infty }{\frac {(a)_{k}(b)_{k}}{(c)_{k}}}{\frac {z^{k}}{k!}}=1+{\frac {ab}{c}}{\frac {z}{1!}}+{\frac {a(a+1)b(b+1)}{c(c+1)}}{\frac {z^{2}}{2!}}+\cdots .
\end{align}
Here $(q)_k$ is the (rising) Pochhammer symbol, which is defined by:
\begin{align}
(q)_{k}=\begin{cases}1&k=0\\q(q+1)\cdots (q+k-1)&k>0\end{cases}
\end{align}
Based on hypergeometric function, we introduce the CDF
\begin{align}
F(x; n)=\frac{n\sinh (x)}{2^n} \,\frac{_2F_1\left(\frac{1}{2},\frac{n+1}{2};\frac{3}{2};-\sinh
^2(x)\right)}{ _2F_1\left(\frac{n}{2},n;\frac{n}{2}+1;-1\right)}, \forall x\in [0, +\infty) \label{eq:hgDis}
\end{align}
where $n>0$ is a real number. It can be verified that $F(x;n)$ is a monolithic non-decreasing function and it satisfies $0\leq F(x; n) \leq 1$. For integer $n \in \{1,2,3\}$, the expressions of $F(x;n)$ are relatively simple:
\begin{align}
F(x;1)&=\frac{2 }{\pi }\arctan(\sinh (x)),\notag\\
F(x;2)&=\tanh (x), \notag\\
F(x;3)&=\frac{{\sech}^2(x)}{\pi}\left(2 \sinh (x)+\arctan(\sinh (x))+\cosh (2 x) \arctan(\sinh(x))\right)\notag
\end{align}
The derivative of this CDF has a very simple expression
\begin{align}
\frac{d}{\ud x }\Big(\sinh (x) \,_2F_1\big(\frac{1}{2},\frac{n+1}{2};\frac{3}{2};-\sinh ^2(x)\big)\Big)=\sech^n(x)
\end{align}
Therefore, the PDF is computed as
\begin{align}
f(x; n)=\frac{n}{2^n} \,\frac{\sech^n(x)}{ _2F_1\left(\frac{n}{2},n;\frac{n}{2}+1;-1\right)}
\end{align}
For integer $n\in \{1,2,3,4,5,6\}$, the PDF has the following expressions:
\begin{align}
f(x;1)=\frac{2}{\pi }{\sech}(x),f(x;2)={\sech}^2(x),f(x;3)=\frac{4}{\pi} {\sech}^3(x)\notag\\
f(x;4)=\frac{3}{2} {\sech}^4(x),f(x;5)=\frac{16}{3 \pi } {\sech}^5(x),f(x;6)=\frac{15}{8}{\sech}^6(x).\notag
\end{align}
The PDF and CDF of this continuous hypergeometric distribution are depicted in Figure \ref{fig:hfp}, which shows that the larger value of $n$ is, the steeper is the graph of the CDF. Compared to the Gudermannian distribution, the continuous hypergeometric distribution generalizes the order to any positive real number.

\subsection{Rapid decay distribution}
In statistical mechanics, the canonical partition function $Z$ of a canonical ensemble is given by
\begin{align}
Z=\exp(-\frac{\phi}{kT}),
\end{align}
where $\phi$ is the Helmholtz free energy, $k$ is the Boltzmann constant, and $T$ is the thermodynamic temperature \cite{gibbs1902elementary}. The function $Z$ is a rapidly decaying exponential function, frequently encountered in physics and applied mathematics. It is particularly useful in contexts where quantities transition smoothly between two regimes, for example, aiding in the construction of smooth step functions and the approximation of Heaviside functions \cite{arfken2011mathematical}.

Based on the canonical partition function, we propose the cumulative distribution function defined as
\begin{align}
F(x)=x\big(1-\exp(-\frac{1}{x})\big) ,\,x\in (0,+\infty).\label{eq:newProb4}
\end{align}
It can be easily verified that $\lim_{x\to 0^+} F(x)=0$, $\lim_{x\to +\infty} F(x)=1$ and $F(x)$ is a non-decreasing continuous function. We refer to the distribution as the rapid decay distribution.
For this function, the probability density function is derived as
\begin{align}
f(x)=F'(x)=1-\exp({-1/x})-\frac{\exp({-1/x})}{x}.\label{eq:pdf4}
\end{align}
\begin{figure}[!htb]
\centering
\includegraphics[width=7cm]{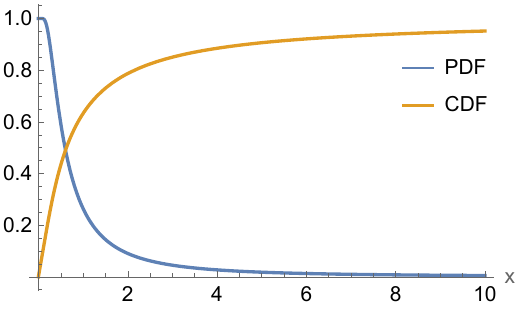}
\caption{Rapid decay distribution: PDF and CDF.}
\label{fig:PDFCDF4}
\end{figure}
The PDF and CDF of this rapid decay distribution are depicted in Figure \ref{fig:PDFCDF4}.
The $m$-th plain moment is
\begin{align}
\mathbb E(x^m)=\int_{0}^{+\infty} f(x) x^m \ud x =m \Gamma(-1-m)\mbox{ if }0<m<1.
\end{align}

\section{More CDF-based damage models}\label{sec:probDamage}

\subsection{Logistic-type damage model based on Logistic distribution}
Replacing the random variable $x$ in Eq.\ref{eq:Logdis} with strain energy density $\phi$ and adding some coefficient
\begin{align}
1-\frac{1}{1+e^{x}}\to \frac{G}{\ell} (1-\frac{2}{1+e^{2\phi\ell/G}})=\frac{G}{\ell} \tanh (\frac{\phi\ell}{G}),\notag
\end{align}
and we arrive at the new fracture energy density defined by
\begin{align}
\psi=\frac{G}{\ell} \tanh (\frac{\phi\ell}{G})
\end{align}
where $\tanh x=\frac{e^{2x}-1}{e^{2x}+1}$ is the hyperbolic tangent function. We designate this innovative model as the logistic-type damage model. 

When $\phi\ell/G$ is small, the Taylor series at $\phi=0$,
\begin{align}
\psi=\phi-\frac{\ell^2}{3 G^2}\phi^3+O(\phi^4),
\end{align}
indicating that the model represents linear elastic behavior when no damage is present.

\begin{figure}[!htb]
\centering
\includegraphics[width=8cm]{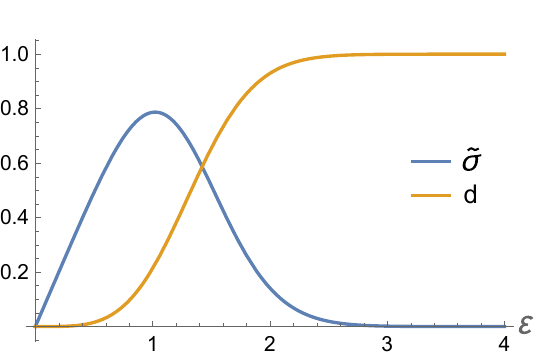}
\caption{Damage model based on half-normal distribution in 1D with $G/\ell=1, k=1$: damage $d$ and effective stress $\tilde{\sigma}=\sech^2(\frac{\phi\ell}{G})\sigma $.}
\label{fig:DamageCrack3}
\end{figure}

The variation of $\psi$ is
\begin{align}
\delta \psi=\frac{4 e^{{2\phi \ell}/{G}}}{(1+e^{{2\phi\ell }/{G}})^2} \delta \phi=\sech^2(\frac{\phi\ell}{G}) \delta \phi
\end{align}
where $\mbox{sech } x=\frac {2e^{x}}{e^{2x}+1}$ is the Hyperbolic secant function and the relation $\frac{d}{d x } \tanh x=1-\tanh^2 x=\sech^2 x$ is used.
We employ the energy functional $E  =\int_\Omega \frac{G}{\ell} \tanh (\frac{\phi^+\ell}{G})+\phi^--\bm b\cdot \bm u \, \ud V$. Through variational derivation of the energy functional, the governing equations are
\begin{align}
\nabla \cdot \big(\sech^2(\frac{\phi\ell}{G})\bm \sigma^++\bm \sigma^-\big)+\bm b=\bm 0
\end{align}
For the convenience of post-processing, we introduce the damage variable defined by
\begin{align}
d=1-\frac{4 e^{2 \phi\ell/G }}{(1+e^{2 \phi\ell/G })^2}=1-\sech^2(\frac{\phi\ell}{G}) ={\tanh}^2(\frac{\phi\ell}{G})
\end{align}
It can be easily verified that in 1D setting
\begin{align}
\int_{0}^{+\infty} \sech^2(\frac{\phi\ell}{G}) \, \sigma \,\ud \varepsilon=\frac{G}{\ell},
\end{align}
which confirms that the envelope of the load-displacement is a constant.
Using relation $\frac{\partial \sigma}{\partial \varepsilon}=0$, the maximal stress is $\bm \sigma_{max}\approx 0.787092 \sqrt{\frac{G k}{\ell}}$ at strain value $\varepsilon\approx 1.02158\sqrt{\frac{G}{k \ell}}$. The corresponding damage is $d\approx 0.229535$.
In order to determine the value of $\ell$, we employ the maximal stress determined from experiment. Similar to the routine in Eq.\ref{eq:Gellv}, the value of $G/\ell$ in current model can be solved from $\bm \sigma_{max}\approx 0.787092 \sqrt{\frac{G k}{\ell}}$ as
\begin{align}
\frac{G}{\ell}=0.619514 \frac{\sigma_{max}^2}{k}. \label{eq:Gellv2}
\end{align}
where $E$ is the elastic modulus.
In the setting of 1D with assumption $G/\ell=1, k=1$ and energy density $\phi=\frac 12 k\varepsilon^2$, we plot the graphs of damage function and effective stress in Figure \ref{fig:DamageCrack3}, which shows non-decreasing property of damage and the stress soften due to damage.
\subsection{Damage model based on half-normal distribution}
Replacing the random variable $x$ in Eq.\ref{eq:errordis} with strain energy density $\phi$ and adding some coefficient
\begin{align}
\erf(x)=\frac{2}{\sqrt{\pi}}\int_0^x e^{-t^2} d t \to \frac{G}{\ell} \erf(\frac{\sqrt{\pi} \ell}{2 G}\phi),\notag
\end{align}
and we arrive at the new fracture energy density defined by
\begin{align}
\psi=\frac{G}{\ell} \erf(\frac{\sqrt{\pi} \ell}{2 G}\phi)
\end{align}
This damage model is named the error-function based damage model or half-normal-type damage model. 

\begin{figure}[!htb]
\centering
\includegraphics[width=8cm]{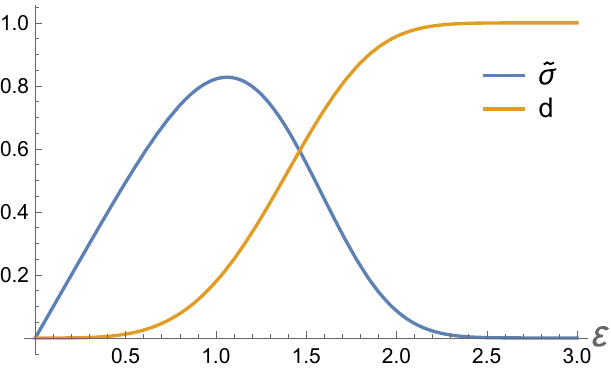}
\caption{Damage model based on the distribution of error function in 1D with $G/\ell=1, k=1$: damage $d$ and effective stress $\tilde{\sigma}=e^{-\frac{\pi \ell^2}{4 G^2} \phi^2}\sigma $.}
\label{fig:DamageCrack5}
\end{figure}

The Taylor series of $\psi$ at point $\phi=0$ is
\begin{align}
\psi=\phi-\frac{\pi \ell^2}{12 G^2} \phi^3+O(\phi^5), 
\end{align}
which highlights the model's ability to represent linear elasticity in damage-free scenarios. 
The variation of $\psi$ is
\begin{align}
\delta \psi=e^{-\frac{\pi \ell^2}{4 G^2} \phi^2} \delta \phi
\end{align}
We employ the energy functional $E  =\int_\Omega\frac{G}{\ell} \erf(\frac{\sqrt{\pi} \ell}{2 G}\phi^+)+\phi^--\bm b\cdot \bm u \, \ud V$. Through variational derivation of the energy functional, the governing equations are
\begin{align}
\nabla \cdot \big(e^{-\frac{\pi \ell^2}{4 G^2} \phi^2}\bm \sigma^++\bm \sigma^-\big)+\bm b=\bm 0
\end{align}
The damage variable defined by
\begin{align}
d=1-e^{-\frac{\pi \ell^2}{4 G^2} \phi^2}
\end{align}
It can be easily verified that in 1D setting
\begin{align}
\int_{0}^{+\infty} e^{-\frac{\pi \ell^2}{4 G^2} \phi^2} \, \sigma \,\ud \varepsilon=\frac{G}{\ell}
\end{align}
Using relation $\frac{\partial \sigma}{\partial \varepsilon}=0$, the maximal stress is $\bm \sigma_{max}=\frac{\sqrt{2}}{\sqrt[4]{e \pi }} \sqrt{\frac{G k}{\ell}}$ at strain value $\varepsilon=\frac{\sqrt{2}}{\sqrt[4]{\pi }} \sqrt{\frac{G}{k \ell}}$. The corresponding damage is $d=1-\frac{1}{\sqrt[4]{e}}$.
In order to determine the value of $\ell$, we employ the maximal stress determined from experiment. Similar to the routine in Eq.\ref{eq:Gellv}, the value of $G/\ell$ in current model can be solved from $\bm \sigma_{max}=\frac{\sqrt{2}}{\sqrt[4]{e \pi }} \sqrt{\frac{G k}{\ell}}$ as
\begin{align}
\frac{G}{\ell}=\frac{2}{\sqrt{e \pi }} \frac{\sigma_{max}^2}{k}. \label{eq:Gellv3}
\end{align}
In the setting of 1D with assumption $G/\ell=1, k=1$ and energy density $\phi=\frac 12 k\varepsilon^2$, we plot the graphs of damage function and effective stress in Figure \ref{fig:DamageCrack5}, which shows non-decreasing property of damage and the stress soften due to damage.
\subsection{Gudermannian-type damage model}
Replacing the random variable $x$ in Eq.\ref{eq:gddis} with strain energy density $\phi$ and adding some coefficient
\begin{align}
\int_{0}^x \sech(\frac{\pi}{2}t)\, dt=\frac{4}{\pi} \arctan(\tanh(\frac{\pi}{4}x))\to \frac{4G}{\pi \ell} \arctan(\tanh(\frac{\pi \ell}{4 G}\phi)),\notag
\end{align}
and we arrive at the new fracture energy density defined by
\begin{align}
\psi=\frac{4G}{\pi \ell} \arctan(\tanh(\frac{\pi \ell}{4 G}\phi))
\end{align}
We term this model the Gudermannian-type damage model. 

\begin{figure}[!htb]
\centering
\includegraphics[width=8cm]{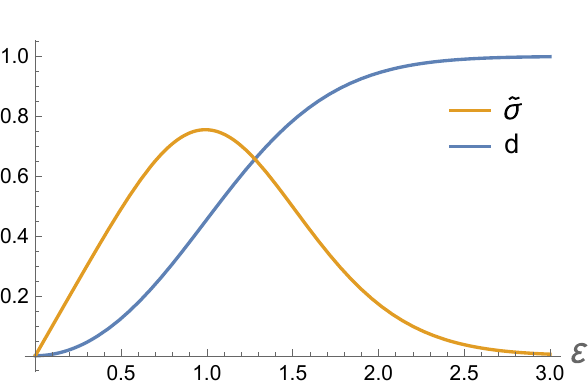}
\caption{Damage model based on Gudermannian distribution in 1D with $G/\ell=1, k=1$: damage $d$ and effective stress $\tilde{\sigma}=\sech(\frac{\phi\ell}{G})\sigma $.}
\label{fig:DamageCrack4}
\end{figure}

The Taylor series of $\psi$ at point $\phi=0$ is
\begin{align}
\psi=\phi-\frac{\pi^2 \ell^2}{24 G^2} \phi^3+O(\phi^5),
\end{align}
which indicates a good linear approximation of $\phi$ at small deformation when the damage is negligible.
The variation of $\psi$ is
\begin{align}
\delta \psi=\sech(\frac{\pi \ell}{2 G} \phi) \delta \phi
\end{align}
We employ the energy functional $E  =\int_\Omega \frac{4G}{\pi \ell} \arctan(\tanh(\frac{\pi \ell}{4 G}\phi^+))+\phi^--\bm b\cdot \bm u \, \ud V$. Through variational derivation of the energy functional, the governing equations are
\begin{align}
\nabla \cdot \big(\sech(\frac{\pi \ell}{2 G} \phi)\bm \sigma^++\bm \sigma^-\big)+\bm b=\bm 0
\end{align}
The damage variable defined by
\begin{align}
d=1-\sech(\frac{\pi \ell}{2 G} \phi)
\end{align}
It can be easily verified that in 1D setting
\begin{align}
\int_{0}^{+\infty} \sech(\frac{\pi \ell}{2 G} \phi) \, \sigma \,\ud \varepsilon=\frac{G}{\ell}
\end{align}
Using relation $\frac{\partial \sigma}{\partial \varepsilon}=0$, the maximal stress is $\bm \sigma_{max}\approx 0.755039\sqrt{\frac{G k}{\ell}}$ at strain value $\varepsilon\approx 0.991243\sqrt{\frac{G}{k \ell}}$. The corresponding damage is $d\approx 0.23829$.
In order to determine the value of $\ell$, we employ the maximal stress determined from experiment. Similar to the routine in Eq.\ref{eq:Gellv}, the value of $G/\ell$ in current model can be solved from $\bm \sigma_{max}=0.755039\sqrt{\frac{G k}{\ell}}$ as
\begin{align}
\frac{G}{\ell}=0.570084\frac{\sigma_{max}^2}{k}. \label{eq:Gellv4}
\end{align}
In the setting of 1D with assumption $G/\ell=1, k=1$ and energy density $\phi=\frac 12 k\varepsilon^2$, we plot the graphs of damage function and effective stress in Figure \ref{fig:DamageCrack4}, which shows non-decreasing property of damage and the stress soften due to damage.

\subsection{Hypergeometric-type damage model}
Based on continuous Hypergeometric distribution in Eq.\ref{eq:hgDis}, do the replacement of $x$ with $\phi \ell/G$
\begin{align}
&\frac{n\sinh (x)}{2^n} \,\frac{_2F_1\left(\frac{1}{2},\frac{n+1}{2};\frac{3}{2};-\sinh
^2(x)\right)}{ _2F_1\left(\frac{n}{2},n;\frac{n}{2}+1;-1\right)}\quad\to\notag\\
&\quad \frac{G}{\ell F_{12}^{(n)}}\sinh \big(\frac{\ell F_{12}^{(n)}}{G}\phi \big)\,_2F_1\Big(\frac{1}{2},\frac{n+1}{2};\frac{3}{2};-\sinh^2\big(\frac{\ell F_{12}^{(n)}}{G}\phi\big)\Big)\notag
\end{align}
where $F_{12}^{(n)}=\frac{2^n}{n}\, _2F_1\left(\frac{n}{2},n;\frac{n}{2}+1;-1\right)$. Some special cases of $F_{12}^{(n)}$ are
\begin{align}
F_{12}^{(0.5)}=4 \sqrt{\frac{2}{\pi }} \Gamma \left(\frac{5}{4}\right)^2,F_{12}^{(1)}=\frac{\pi }{2},F_{12}^{(2)}=1,\notag\\
F_{12}^{(3)}=\frac{\pi }{4},F_{12}^{(4)}=\frac{2}{3},F_{12}^{(5)}=\frac{3 \pi }{16},F_{12}^{(6)}=\frac{8}{15}\notag
\end{align}
We define the fracture energy density based on $\phi^+$, the ``positive'' part of the energy density $\phi$:
\begin{align}
\psi= \frac{G}{\ell F_{12}^{(n)}}\sinh \big(\frac{\ell F_{12}^{(n)}}{G}\phi^+ \big)\,_2F_1\Big(\frac{1}{2},\frac{n+1}{2};\frac{3}{2};-\sinh^2\big(\frac{\ell F_{12}^{(n)}}{G}\phi^+\big)\Big)\label{eq:hfdm}
\end{align}
where $n>0$ is a real number. It can be verified that $0\leq \psi\leq \frac{G}{\ell}$.
The Taylor expansion of $\psi$ at $\phi^+=0$ is
\begin{align}
\psi=\phi^+ -\frac{n}{6} \left(\frac{F_{12}^{(n)} \ell}{ G}\right)^2\phi^{+3}+O\left(\phi ^{+5}\right) \label{eq:hfdmT}
\end{align}
When $\phi$ is small, $\psi\approx \phi^+$, which confirms that the state without damage is well approximated.
The variation of $\psi$ has a very compact form
\begin{align}
\delta \psi=\sech^n \big(\frac{\ell F_{12}^{(n)}}{G}\phi^+\big) \delta \phi^+
\end{align}
We employ the energy functional $E  =\int_\Omega \psi+\phi^--\bm b\cdot \bm u \, \ud V$. Through variational derivation of the energy functional, the governing equations are
\begin{align}
\nabla \cdot \big(\sech^n \big(\frac{\ell F_{12}^{(n)}}{G}\phi^+\big)\bm \sigma^++\bm \sigma^-\big)+\bm b=\bm 0
\end{align}
The damage variable is defined by
\begin{align}
d=1-\sech^n \big(\frac{\ell F_{12}^{(n)}}{G}\phi^+\big)\label{eq:hfmd1}
\end{align}
This model is identified as the hypergeometric-type damage model.

\begin{figure}[htp]
\centering
\subfigure[]{
\includegraphics[width=0.45\textwidth]{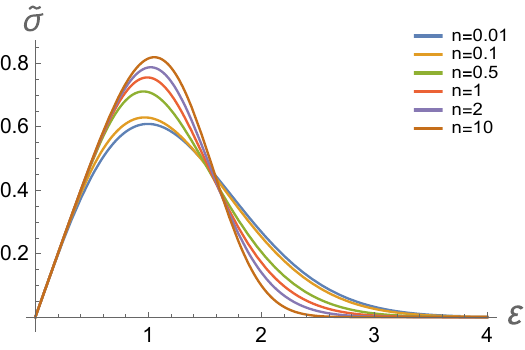}
\label{fig:hfsigma}}
\subfigure[]{
\includegraphics[width=0.45\textwidth]{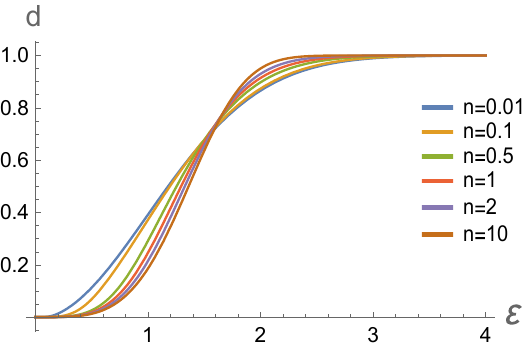}
\label{fig:hfdamage}}
\caption{Hypergeometric-type damage model by Eq.\ref{eq:hfdm} in 1D: (a) the effective stress versus strain; (b) the damage versus strain based on Eq.\ref{eq:hfmd1}.}
\label{fig:hfdm}
\end{figure}

In the setting of 1D with assumption $G/\ell=1, k=1$ and energy density $\phi=\frac 12 k\varepsilon^2$, the evolution of effective stress and damage of the hypergeometric damage model is illustrated in Figure \ref{fig:hfdm}.
When $n=1$, the Hypergeometric-type damage model reduces to the Gudermannian damage model. When $n=2$, the Hypergeometric-type damage model degenerates to the Logistic damage model. In general, the Hypergeometric-type damage model generalizes these two models to any positive real order.

\subsection{Damage model based on piecewise distribution}
\begin{figure}[!htb]
\centering
\subfigure[]{
\includegraphics[width=0.45\textwidth]{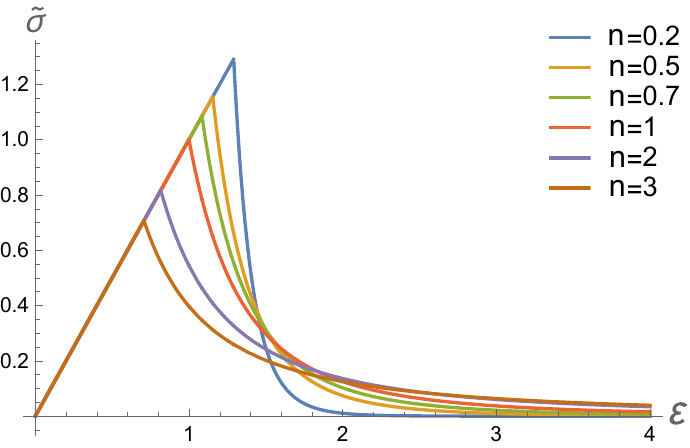}
\label{fig:sigmaEffectiveMs1}}
\subfigure[]{
\includegraphics[width=0.45\textwidth]{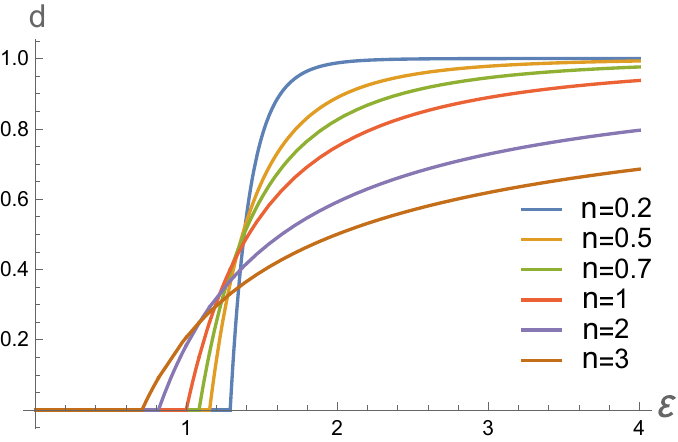}
\label{fig:damageEffectiveM0b}}
\caption{Damage model by Eq.\ref{eq:VDthresholdA} in 1D: (a) the effective stress versus strain; (b) the damage versus strain based on Eq.\ref{eq:vdmd2b}}
\label{fig:EffectiveMs0b}
\end{figure}

Using Eq. \ref{eq:newProb2} and replacing the random variable $x$ with $\phi\ell/G$,e.g.
\begin{align}
\mbox{If}\big(x\leq \frac{1}{n+1}, x, 1-\frac{{n}/(n+1)}{((n+1)x)^{1/n}}\big)\,\to \mbox{If}\big(\phi\frac{\ell}{G}\leq \frac{1}{n+1}, \phi, \frac{G}{\ell}\big( 1-\frac{{n}/(n+1)}{((n+1)\phi\ell/G)^{1/n}}\big)\big)
\end{align}
Accordingly, we define the fracture energy density as
\begin{align}
\psi= \mbox{If}\big(\phi\ell/G\leq \frac{1}{n+1}, \phi, \frac{G}{\ell}\big( 1-\frac{{n}/(n+1)}{((n+1)\phi\ell/G)^{1/n}}\big)\big)
\end{align}
and use the energy functional $E  =\int_\Omega \psi+\phi^--\bm b\cdot \bm u \, \ud V$. Based on variational derivation and integration by parts, the governing equations of elastic solid embedded with damage model are
\begin{align}
\nabla\cdot \Big(\min(1, \big(\frac{G }{(n+1)\phi^+ \ell}\big)^{\frac{n+1}{n}}) \,\bm \sigma^++\bm \sigma^-\Big)+\bm b=\bm 0.\label{eq:VDthresholdA}
\end{align}
In Eq.\ref{eq:VDthresholdA}, the damage variable for the purpose of post-processing is defined as
\begin{align}
d=\max(0, 1-\big(\frac{G }{(n+1)\phi^+ \ell}\big)^{\frac{n+1}{n}})\label{eq:vdmd2b}
\end{align}
The damage model given by Eq.\ref{eq:VDthresholdA} coincides with the variational damage model in Ref \cite{Ren2024vdm2}.
In the current model, the crack length scale remains undetermined. In order to determine the value of $\ell$, we employ the maximal stress determined from experiment. Thus, we solve $G/\ell$ from $\bm \sigma_{max}$ as
\begin{align}
\frac{ G}{\ell}= \frac{n+1}{2}\frac{\sigma_{max}^2}{k}.
\end{align}
In the setting of 1D with assumption $G/\ell=1, k=1$ and energy density $\phi=\frac 12 k\varepsilon^2$, the evolution of effective stress and damage is illustrated in Figure \ref{fig:EffectiveMs0b}
\subsection{Damage model based on rational distribution}
Using Eq. \ref{eq:newProb3} and replacing the random variable $x$ with $\phi\ell/G$,e.g.
\begin{align}
\frac{x^2+x n^2}{(x+n)^2}\,\to \,\frac{\phi\ell/G+n^2}{(\phi\ell/G+n)^2} \,\,\phi
\end{align}
where $0.5\leq n\leq 2$.
Accordingly, we define the fracture energy density as
\begin{align}
\psi=\frac{\phi^+\ell/G+n^2}{(\phi^+\ell/G+n)^2} \,\,\phi^+
\end{align}
The Taylor expansion of $\psi$ at $\phi^+=0$ is
\begin{align}
\psi=\phi^+ -\frac{(2n-1)}{n^2} \frac{\ell}{G}\phi^{+2}+\frac{(3n-2)}{n^3} \frac{\ell^2}{G^2}\phi^{+3}+O\left(\phi ^{+4}\right),
\end{align}
which reflects the model's consistency with linear elasticity for undamaged materials.

\begin{figure}[!htb]
\centering
\subfigure[]{
\includegraphics[width=0.45\textwidth]{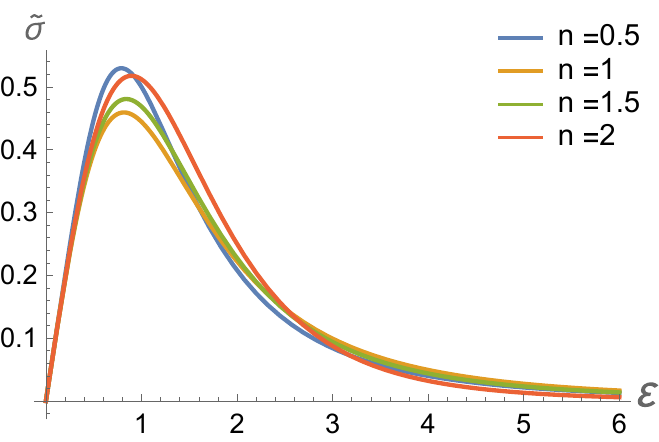}
\label{fig:sigmaCase10}}
\subfigure[]{
\includegraphics[width=0.45\textwidth]{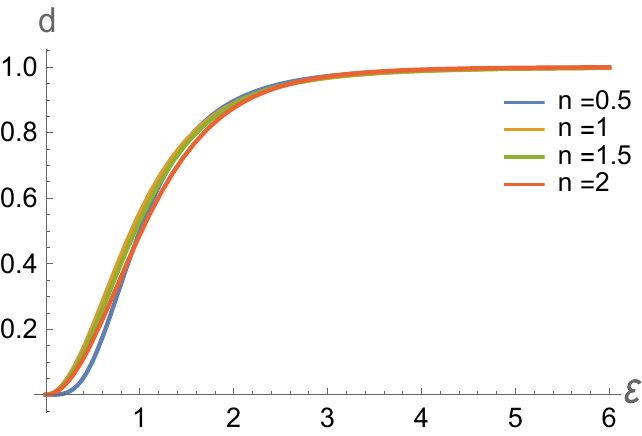}
\label{fig:damageCase10}}
\caption{Damage model by Eq.\ref{eq:VDthresholdA} in 1D: (a) the effective stress versus strain; (b) the damage versus strain based on Eq.\ref{eq:vdmd2bR}}
\label{fig:EffectiveMs0b10}
\end{figure}

The variation of $\psi$ is
\begin{align}
\delta \psi=\frac{{(2n-n^2) \phi^+\ell }/{G}+n^3}{({\phi^+ \ell}/{G}+n)^3} \,\,\delta \phi^+
\end{align}
We use the energy functional $E  =\int_\Omega \psi+\phi^--\bm b\cdot \bm u \, \ud V$. Based on variational derivation and integration by parts, the governing equations of elastic solid embedded with damage model are
\begin{align}
\nabla\cdot \Big(\frac{{(2n-n^2) \phi^+\ell }/{G}+n^3}{({\phi^+ \ell}/{G}+n)^3} \,\bm \sigma^++\bm \sigma^-\Big)+\bm b=\bm 0.\label{eq:VDMration}
\end{align}
In Eq.\ref{eq:VDthresholdA}, the damage variable for the purpose of post-processing is defined as
\begin{align}
d=1-\frac{{(2n-n^2) \phi^+\ell }/{G}+n^3}{({\phi^+ \ell}/{G}+n)^3}\label{eq:vdmd2bR}
\end{align}
In the setting of 1D with assumption $G/\ell=1, k=1$ and energy density $\phi=\frac 12 k\varepsilon^2$, the evolution of effective stress and damage is illustrated in Figure \ref{fig:EffectiveMs0b10}

\subsection{Damage model based on rapid decay distribution}
Using Eq. \ref{eq:newProb4} and replacing the random variable $x$ with $\phi\ell/G$,e.g.
\begin{align}
x(1-\exp(-1/x))\,\to \, \phi(1-\exp(-\frac{G}{\phi\ell}))
\end{align}
Accordingly, we define the fracture energy density as
\begin{align}
\psi=(1-\exp(-\frac{G}{\ell\phi^+})) \phi^+
\end{align}
Interestingly, this fracture energy functional is quite close to the hyperelasticity model with a softening effect in Ref \cite{volokh2007hyperelasticity,tran2024new}. When $\phi$ is small, $\psi\approx \phi^+$ as $\phi^+ \exp(-\frac{G}{\ell \phi^+})\approx 0$, which confirms that the state without damage is well approximated.
The variation of $\psi$ is
\begin{align}
\delta \psi=\big(1-(1+\frac{G}{\ell\phi^+})\exp(-\frac{G}{\ell \phi^+})\big) \,\,\delta \phi^+
\end{align}
We use the energy functional $E  =\int_\Omega \psi+\phi^--\bm b\cdot \bm u \, \ud V$. Based on variational derivation and integration by parts, the governing equations of elastic solid embedded with current damage model are
\begin{align}
\nabla\cdot \Big(\big(1-(1+\frac{G}{\ell\phi^+})\exp(-\frac{G}{\ell \phi^+})\big) \,\bm \sigma^++\bm \sigma^-\Big)+\bm b=\bm 0.\label{eq:VDMrapidDecay}
\end{align}
In above equation, the damage variable for the purpose of post-processing is defined as
\begin{align}
d=(1+\frac{G}{\ell\phi^+})\exp(-\frac{G}{\ell \phi^+}).
\label{eq:VDMrapidDecayDmg}
\end{align}

\bibliographystyle{unsrt}
\bibliography{VDMprobability}
\end{document}